\documentclass[a4paper,
10pt
]{article}

\usepackage[english]{babel}

\usepackage[a4paper,left=3cm,right=3cm,top=2.5cm,bottom=2.5cm]{geometry}

\usepackage{hyperref}
\usepackage{bm}
\usepackage{authblk}
\usepackage{siunitx}
\usepackage{graphicx}
\usepackage{amsmath}
\graphicspath{ {./images/} }
\usepackage[table,xcdraw]{xcolor}
\usepackage{hhline}
\definecolor{lightgray}{gray}{0.9}
\usepackage{tikz}
\usepackage{subcaption}
\usetikzlibrary{positioning}
\usepackage[draft,inline,marginclue]{fixme}
\usepackage[numbers,sort&compress]{natbib}
\usepackage{subcaption}

\usepackage{float} 
\usepackage{rotating} 

\usepackage{booktabs}
\usepackage[table,xcdraw]{xcolor}
\usepackage{array}
\newcolumntype{P}[1]{>{\raggedright\arraybackslash}p{#1}}

\FXRegisterAuthor{gs}{adt}{\color{red}GS}

\title{Stochastic Parameter Prediction in Cardiovascular Problems}

\author[1]{Kabir ~Bakhshaei\footnote{kbakhsha@sissa.it}}
\author[1]{Sajad ~Salavatidezfouli\footnote{ssalavat@sissa.it }}
\author[2]{Giovanni ~Stabile\footnote{giovanni.stabile@santannapisa.it}}
\author[1]{Gianluigi ~Rozza\footnote{grozza@sissa.it: Corresponding author}}

\affil[1]{Mathematics Area, mathLab, SISSA, Trieste, Italy}
\affil[2]{The Biorobotics Institute, Sant'Anna School of Advanced Studies, Pisa, Italy}

\date{\today} 

\begin{document}
\maketitle
\listoffixmes
\begin{abstract}

Patient-specific modeling of cardiovascular flows with high-fidelity is challenging due to its dependence on accurately estimated velocity boundary profiles, which are essential for precise simulations and directly influence wall shear stress calculations—key in predicting cardiovascular diseases like atherosclerosis. This data, often derived from in-vivo modalities like 4D flow MRI, suffers from low resolution and noise. To address this, we employ a stochastic data assimilation technique that integrates computational fluid dynamics with an advanced Ensemble-based Kalman filter enhancing model accuracy while accounting for uncertainties. Our approach sequentially collects velocity data over time within the vascular model, enabling real-time refinement of unknown boundary estimations. The mathematical model uses the incompressible Navier–Stokes equation to simulate aortic blood flow. We consider unknown boundaries as constant, time-dependent, and time-space-dependent in two- and three-dimensional models. In our 2-dimensional model, relative errors were as low as 0.996\% for constant boundaries and up to 2.63\% and 2.61\% for time-dependent boundaries and time-space-dependent boundaries, respectively, over an observation span of two-time steps. For the 3-dimensional patient-specific model, the relative error was 7.37\% for space-time-dependent boundaries. By refining the velocity boundary profile, our method improves wall shear stress predictions, enhancing the accuracy and reliability of models specific to individual cardiovascular patients. These advancements could contribute to better diagnosis and treatment of cardiovascular diseases.\\

Keywords: Stochastic Data Assimilation, Ensemble Kalman Filter, Cardiovascular Flows, Bayesian Inversion, Uncertainty Quantification, Computational Hemodynamics.
\end{abstract}

\section{Introduction}\label{sec:Intro}

To predict blood flow velocities and vessel geometries, in vivo methods in a clinical setting typically utilize non-invasive imaging methods like ultrasound or magnetic resonance imaging (MRI) \citep{agbalessi2023tracking}. However, it is challenging to directly measure in vivo wall shear stress (WSS) employing conventional clinical methods, despite its critical role in predicting cardiovascular diseases such as aneurysm and atherosclerosis \citep{urschel2021investigation}. Consequently, in-vivo tests alone cannot provide the predictive capabilities that can be achieved through the use of simulations of intricate cardiovascular systems \citep{peirlinck2021precision}.\\

Computational cardiology has made enormous advances in the last ten years, despite being a multi-disciplinary field \citep{kashtanova2023simultaneous, bai2023diagnosis} . Specifically, researchers have been able to tackle important problems associated with the limitations of clinical methods thanks to new developments in numerical analysis, computational resources, and patient-specific models called digital twins \citep{de2022evidence}. Because blood flow patterns are known to have a significant impact on potentially life-threatening cardiovascular diseases such as aneurysms\citep{salavatidezfouli2024effect}, stenosis\citep{nadeem2023numerical}, and aortic dissection \citep{wu2023vitro} hemodynamics simulations have been useful for assessing these conditions.\\

A variety of hemodynamic models have been developed to analyze both normal and pathological physiologies, ranging from electric analogues of the circulation, which are lumped models, to 1-dimensional Euler equations and 3-dimensional Navier-Stokes models, and finally, geometrical multi-scale models \citep{pant2017inverse}. The model parameters, including blood properties and inflow/outflow boundary conditions, must be supplied to all the methods that have been mentioned so far \citep{armour2021influence}. A medical imaging modality called Phase Contrast Magnetic Resonance Imaging (PC-MRI) \citep{stankovic20144d} was developed to visualize and quantify the blood flow in the vasculature. The most typical method for obtaining boundary parameters is based on this research \citep{martin20242115}. This imaging method is referred to as 4D PC-MRI or 4D-flow MRI in the literature \citep{zhuang2021role}. It encompasses time-resolved slice images that span an acquisition volume and employ flow encoding. It allows for the non-invasive examination of blood hemodynamics in the vessel by providing both anatomical (magnitude images) and functional (phase images) information \citep{morgan20214d}. \\

Nevertheless, the process of generating the velocity profile in a specific section of the vessel using 4D-flow analysis entails significant pre-processing \citep{agbalessi2023tracking}. Specifically, it may be necessary to extract the properties of the aorta at each time interval of the cardiac cycle \citep{tahoces2019automatic}. Additionally, 4D MRI data is inherently noisy and uncertain, partly due to physiological motions, such as heartbeats and breathing, which introduce natural variability and measurement errors during data acquisition. These factors necessitate dedicated pre-processing steps to reduce noise and increase reliability \citep{de2020data}. Consequently, accurately determining the model parameters and boundary profile of the vascular system remains uncertain \citep{arthurs2020flexible}. For the current study, we employ Data assimilation(DA) techniques to overcome these limitations by integrating distributed in-vivo or in-vitro data \citep{bertoglio2014identification}. This approach mitigates the noise impact and enhances simulation accuracy, thereby providing a more precise representation of cardiovascular dynamics. 
\\

The DA method has become increasingly popular in recent years due to its ability to integrate continuously noisy measurements into a mathematical model \citep{carrassi2018data}. This method adjusts model parameters and states in real-time while addressing and reducing both model and parameter uncertainty in predictions by continuously being updated. This approach also significantly improves the accuracy and reliability of forecasts and simulations in fields like meteorology\citep{wang2024data}, oceanography\citep{valmassoi2023current}, and hydrology \citep{shao2024wrf}. The utilization of observations in conjunction with model estimations in this method effectively minimizes uncertainties and enhances the efficacy of decision-making processes \citep{de2021data}.  It also has been used in various cardiovascular applications \citep{lim2023data}, including the prediction of Windkessel parameters \citep{fevola2023vector}, the determination of closed-loop lumped model parameters for single-ventricle physiology \citep{bjordalsbakke2022parameter}, the assessment of artery wall stiffness parameters \citep{imperiale2021sequential}, and the estimation of tissue support parameters for fluid-structure interaction models \citep{shidhore2022estimating}.\\

Numerous prior studies have addressed the inverse problem, specifically concentrating on predicting the parameters that establish the boundary conditions of the model \citep{arbia2015pulmonary, perdikaris2016model, fevola2022boundary, SmithThesis}. A novel Bayesian approach has been recently devised for estimating parameters in cardiovascular models, using methods based on Markov chain Monte Carlo (MCMC) to estimate mean or maximum values, as described in reference \citep{Schiavazzi2017Patient}. Sequential DA techniques, such as Kalman Filter(KF) variants are methods that stand out from others due to their flexibility and computational efficiency in estimating parameters in systems. They achieve this by dynamically combining observational data with predictive models, resulting in improved accuracy as time progresses \citep{Canuto2020AnEK, Walia2023PPGbasedWM}.\\ 


This paper presents a technique for estimating parameters in patient-specific cardiovascular models, which is then evaluated through testing. The primary method used for estimation is an advanced variant of the ensemble Kalman filter (EnKF), initially created to assimilate data in dynamic models for numerical weather prediction\citep{Fang2017EnsemblebasedSI}. The main benefit of the EnKF lies in its suitability for nonlinear dynamical systems. Although other methods are available for estimating parameters in cardiovascular models \citep{Canuto2020AnEK}, the literature does not currently report any similar work for high-fidelity vascular models. This absence is mainly due to challenges in implementing the detailed models, managing their complex dynamics, and meeting the computational requirements. Most existing studies employ simplified or linearized models rather than directly solving the full Navier-Stokes equations, often using approximations in either model geometry or dynamics to reduce computational complexity. By contrast, our approach seeks to capture these dynamics with a high degree of fidelity to patient-specific conditions, providing a more robust and accurate representation of cardiovascular behavior. The key contributions and findings of this paper are outlined as follows.
\begin{itemize}
    \item Integration of an advanced ensemble Kalman filter extension for simultaneously estimating the system's states and inputs in cardiovascular inverse problem with stochastic Bayesian approach.
    \item Investigation of various distributions of unknown velocity boundary conditions, including constant, time-dependent, and time-space-dependent cases, and assessment of the proposed ensemble model on both a 2-dimensional idealized model and a 3-dimensional patient-specific model to solve the filter’s forward problem.
    \item Demonstration of the robustness of the data assimilation method, highlighting the proposed approach's strength in delivering accurate predictions using a low-fidelity CFD solver for the forward problem. 
    \item Finally, assessment of the proposed EnKF’s ability to reconstruct the flow field and parameter prediction across different flow regimes, with a specific focus on its performance during turbulent intervals, such as at systole peak.
\end{itemize}

The paper follows the following structure: Section \hyperref[sec:EnSISSF]{2} discusses the data assimilation process, including the initial setup, prediction, and update stages. Section \hyperref[sec:Dynamic Model]{3}, applies this approach to cardiovascular flow, detailing the mathematical model, turbulence considerations, boundary conditions, numerical solutions, and the generation of synthetic measurement data for data assimilation. Section \hyperref[sec:resVessel]{4} presents the results of the 2D simulations, exploring scenarios with constant parameters, time-dependent parameters, and time-space-dependent parameters. Section \hyperref[3DTimeSpace]{5} extends these findings to a 3D simulation for the time-space-dependent parameter. Finally, Section \hyperref[sec:conc]{6} provides the conclusion of our study.

\begin{figure}[H]
    \centering
    \includegraphics[width=1\linewidth]{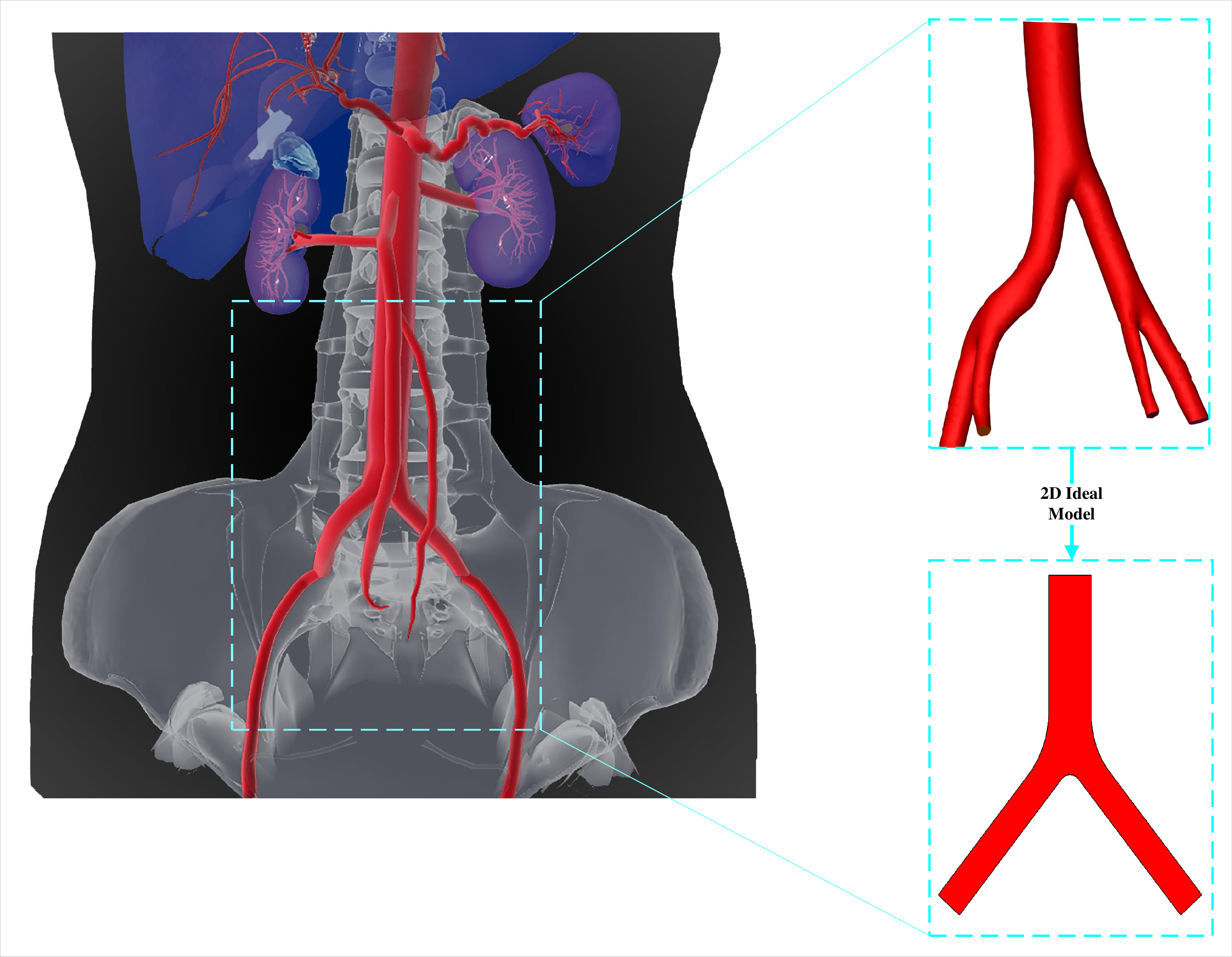}
    \caption{Abdominal aorta with 2D ideal and 3D patient-specific models}
    \label{TransparentAbdominalAorta}
\end{figure}

\section{Data Assimilation}\label{sec:EnSISSF}
In this research, we employ an advanced version of the EnKF, known as the Ensemble-based Simultaneous Input and State Filtering (EnSISF) with direct feedthrough (EnSISF-wDF)\citep{Fang2017-tx, bakhshaei2024optimized}, as a data assimilation method. This approach is designed not only to estimate the unknown boundary velocity profile as an input parameter but also to predict the velocity and pressure within the system. It does so by resolving the transient Navier-Stokes equations for the forward problem of the method. The general nonlinear dynamic system, which allows direct input-to-output feed-through is given by.
\begin{center}
\begin{flalign}
& \left\{
\begin{aligned}
\bm{\chi}_{k+1}=\bm{\Lambda}\left(\bm{\chi}_{k}, \bm{u}_{k}\right) +\bm{w
}_{k}, \\
\bm{y}_{k}=\bm{\Theta}\left(\bm{\chi}_{k}, \bm{u}_{k}\right) +\bm{v
}_{k}.
\end{aligned}
\right.
\label{eq:SISSFwdf}
\end{flalign}
\end{center}

During each time step $k+1$, $\bm{\chi}_{k+1}$ as the unknown system state vector is computed based on the nonlinear transition function $\bm{\Lambda}$, which depends on unknown input vector $\bm{u}_{k}$, and the prior state state $\bm{\chi}_{k}$. The measurement vector, $\bm{y}_{k}$ denotes the observed system values, with $\bm{\Theta}$ representing either a linear or nonlinear measurement function.\\

The process noise $\bm{w}$ and measurement noise $\bm{v}$ are considered independent, zero-mean and are assumed to follow Gaussian distributions characterized by covariance matrices $\bm{Q}_{k}$ and $\bm{R}_{k}$ respectively. Here, $\bm{w}$ captures 
discrepancies due to unexpected or unmodeled variations not accounted for by the mathematical model.Measurement noise $\bm{v}_{k}$, which accounts for erros in the data, can arise from sensor limitations, environmental factors, or external disturbances. These noises are represented as follows:

\begin{center}
\begin{equation}
\bm{w}_{k} = \mathcal{N}(\bm{0}, \bm{Q}_{k}).
\label{eq:ProcessNoise}
\end{equation}
\end{center}

\begin{center}
\begin{equation}
\bm{v}_{k} = \mathcal{N}(\bm{0}, \bm{R
}_{k}).
\label{eq:MeasurementsNoise}
\end{equation}
\end{center}

At each time step k = 1, 2, ..., this method decouples the prediction and update process, assembles them into ensembles, and combine then iteratively. The following sections outline the prediction-update procedure that this study adopts, along with theoretical foundations \cite{Fang2017-tx}.

\subsection{Initial State Configuration Using Gaussian Priors}\label{subsec:Initial}
The procedure begins by setting Gaussian priors for the initial state of the system $\bm{\chi}_{k=0}$, encompassing variables such as the x-velocity, y-velocity, and pressure as defined in Eq. \ref{eq:GausPriorstate}. Additionally, the initial condition for the unknown boundary velocity in Eq. \ref{eq:GausPriorParameter}, is also assigned a Gaussian prior, characterized by a mean $\bm{\hat{u}}_{j,k=0}$ and a covariance represented by $\bm{\kappa}$. The initial system state $\bm{\chi}_{k=0}$ is assume to adhere to a Gaussian distribution with a mean of $\bm{T_0}$ and covariance $\bm{\sigma}$.
\begin{center}
\begin{equation}
\bm{\chi}_{k=0} = \mathcal{N}(\bm{T_0}, \bm{\sigma}_{T}).
\label{eq:GausPriorstate}
\end{equation}
\end{center}

\begin{center}
\begin{equation}
\bm{u}_{k=0} = \mathcal{N}(\bm{\hat{u}}_{j,k=0}, \bm{\kappa}).
\label{eq:GausPriorParameter}
\end{equation}
\end{center}

\subsection{Predictive Estimation Step}\label{subsec:Predict}
In each time step, the forecast phase uses information from the prior time step to predict the current state. The process model $\bm{\Lambda}$ is applied directly, incorporating unknown parameters specific to each ensemble member. The outcome is combined with the model error, as described in Eq. \ref{eq:ProcessNoise}, producing an ensemble of possible solutions as presented in Eq. \ref{eq:forwardProblem}. Due to the stochastic nature of the boundary conditions, the states are also treated stochastically. The forecasted state $\widehat{\bm{\chi}}_{k \mid k-1}^i$ for each ensemble member $i$ at time $k$ is expressed as:
\begin{center}
\begin{equation}
\widehat{\bm{\chi}}_{k \mid k-1}^i=\bm{\Lambda}\left(\hat{\bm{\psi}}_{k-1 \mid k-1}^i\right)+\bm{w}_{k-1}^i .
\label{eq:forwardProblem}
\end{equation}
\end{center}

Here $\widehat{\bm{\psi}}_{k \mid k-1}$ represenr joint ensemble state from the prior time step and $\bm{w}_{k-1}^i$ denotes process noise for the ensemble number $i$. The joint ensemble $\widehat{\bm{\psi}}_{k \mid k-1}$ is structured as follows: 
\begin{center}
\begin{equation}
\hat{\bm{\psi}}_{k \mid k-1}^{i}=\left[\begin{array}{c}
\hat{\bm{u}}_k^{i} \\
\hat{\bm{\chi}}_{k \mid k-1}^{i}
\end{array}\right].
\label{eq:builtJointEns}
\end{equation}
\end{center}

The ensemble mean $\widehat{\bm{\chi}}_{k \mid k-1}$ is then calculated using the following formula, where $S_n$ represents the total number of the ensemble members. 
\begin{center}
\begin{equation}
\widehat{\bm{\chi}}_{k \mid k-1} = \frac{1}{S_n} \sum_{i=1}^{S_n} \widehat{\bm{\chi}}_{k \mid k-1}^i.
\label{eq:forwardProblemMean}
\end{equation}
\end{center}

If measurements are unavailable at a given time step, a combined ensemble is generated following Eq. \ref{eq:builtJointEns}, and its mean is determined using Eq. \ref{eq:meanjointEnsembleUpdate}.

\subsection{Refinement Update Step}\label{subsec:Update}

To enhance the accuracy of the estimation process, measurements are incorporated iteratively, where each iteration is labeled by $\beta$. For each ensemble member, forecasted values at the 
$\beta$th iteration are generated using:
 
\begin{center}
\begin{equation}
\hat{\bm{y}}_{k \mid k-1}^{i, \beta}=\bm{\Theta}\left(\hat{\bm{\psi}}_{k \mid k-1}^{i, \beta}\right)+\bm{v}_k^{i, \beta}.
\label{eq:measurementFunction}
\end{equation}
\end{center}

In the subsequent update phase, once new measurements are obtained, the observational covariance matrix $\bm{P}_{k \mid k-1}^{\bm{y}, \beta}$ and the cross-covariance matrix $\bm{P}_{k \mid k-1}^{\bm{\psi} \bm{y}, \beta}$, which captures the relationship between the joint variables and observations, are determined. This process is guided by Eqs. \ref{eq:meanjointEnsembleUpdate}, and \ref{eq:meanMeasurementsEnsembleUpdate}, which represent the average of the joint $\hat{\boldsymbol{\psi}}_{k \mid k-1}$ and observation ensembles $\hat{\bm{y}}_{k \mid k-1}^{\beta}$, respectively.
\begin{center}
 \begin{equation}
\hat{\boldsymbol{\psi}}_{k \mid k-1}^{\beta}=\frac{1}{S_n} \sum_{i=1}^{S_n} \hat{\bm{\psi}}_{k \mid k-1}^{i, \beta}.
\label{eq:meanjointEnsembleUpdate}
\end{equation}
\end{center}

\begin{center}
\begin{equation}
\hat{\bm{y}}_{k \mid k-1}^{\beta}=\frac{1}{S_n} \sum_{i=1}^{S_n} \hat{\bm{y}}_{k \mid k-1}^{i, \beta}.
\label{eq:meanMeasurementsEnsembleUpdate}
\end{equation}
\end{center}

 The covariance matrices are given by: 
\begin{center}
\begin{equation}
\bm{P}_{k \mid k-1}^{\bm{y}, \beta}=\frac{1}{S_n} \sum_{i=1}^{S_n} \hat{\bm{y}}_{k \mid k-1}^{i, \beta} \hat{\bm{y}}_{k \mid k-1}^{i, \beta \top}-\hat{\bm{y}}_{k \mid k-1}^{\beta} \hat{\bm{y}}_{k \mid k-1}^{\beta \top}.
\label{eq:covarianceMeasurementsSigmaPoint}
\end{equation}
\end{center}

\begin{center}
\begin{equation}
\bm{P}_{k \mid k-1}^{\bm{\psi} \bm{y}, \beta}=\frac{1}{S_n} \sum_{i=1}^{S_n} \hat{\bm{\psi}}_{k \mid k-1}^{i, \beta} \hat{\bm{y}}_{k \mid k-1}^{i, \beta \top}-\hat{\boldsymbol{\psi}}_{k \mid k-1}^{\beta} \hat{\bm{y}}_{k \mid k-1}^{\beta \top}.
\label{eq:crossCovarianceMeasurementsSigmaPointAndSigmaPoints}
\end{equation}
\end{center}

The ensemble Kalman $\bm{K}_E$ gain shown in Eq. \ref{eq:EnsembleKalmanGain} is an essential part of this approach as it enhances the joint variable estimation accuracy. This is achieved by itegrating observational data into the ensemble of model states and input parameters.
\begin{center}
\begin{equation}
\bm{K}_E = \bm{P}_{k \mid k-1}^{\bm{\psi} \bm{y}, \beta}\left(\bm{P}_{k \mid k-1}^{\bm{y},\beta}\right)^{-1}.
\label{eq:EnsembleKalmanGain}
\end{equation}
\end{center}

Eq. \ref{eq:updateJointEnsemble} demonstrates the process of updating the ensemble of joint variables $\hat{\boldsymbol{\psi}}_{k \mid k}^{i, \beta}$. This update is achieved by adding previous joint ensemble $\hat{\boldsymbol{\psi}}_{k \mid k-1}^{i, \beta}$ with the product of $\bm{K}_E$ as the Kalman gain and the difference between the observed data $\bm{y}_k$ and the predicted measurement \ref{eq:measurementFunction}.
\begin{center}
\begin{equation}
\hat{\boldsymbol{\psi}}_{k \mid k}^{i, \beta}=\hat{\boldsymbol{\psi}}_{k \mid k-1}^{i, \beta}+\bm{K}_E\left(\bm{y}_k-\hat{\bm{y}}_{k \mid k-1}^{i, \beta}\right) .
\label{eq:updateJointEnsemble}
\end{equation}
\end{center}

Finally, the mean values of the updated states and input parameters are obtained using Equ. \ref{eq:meanOfUpdateJointEnsemble}.
\begin{center}
\begin{equation}
\hat{\bm{\psi}}_{k \mid k}^{\beta}=\frac{1}{S_n} \sum_{i=1}^{S_n} \hat{\bm{\psi}}_{k \mid k}^{i, \beta}.
\label{eq:meanOfUpdateJointEnsemble}
\end{equation}
\end{center}

\begin{figure}[H]
    \centering    \includegraphics[width=1\textwidth, keepaspectratio, trim=0 0 0 0, clip]{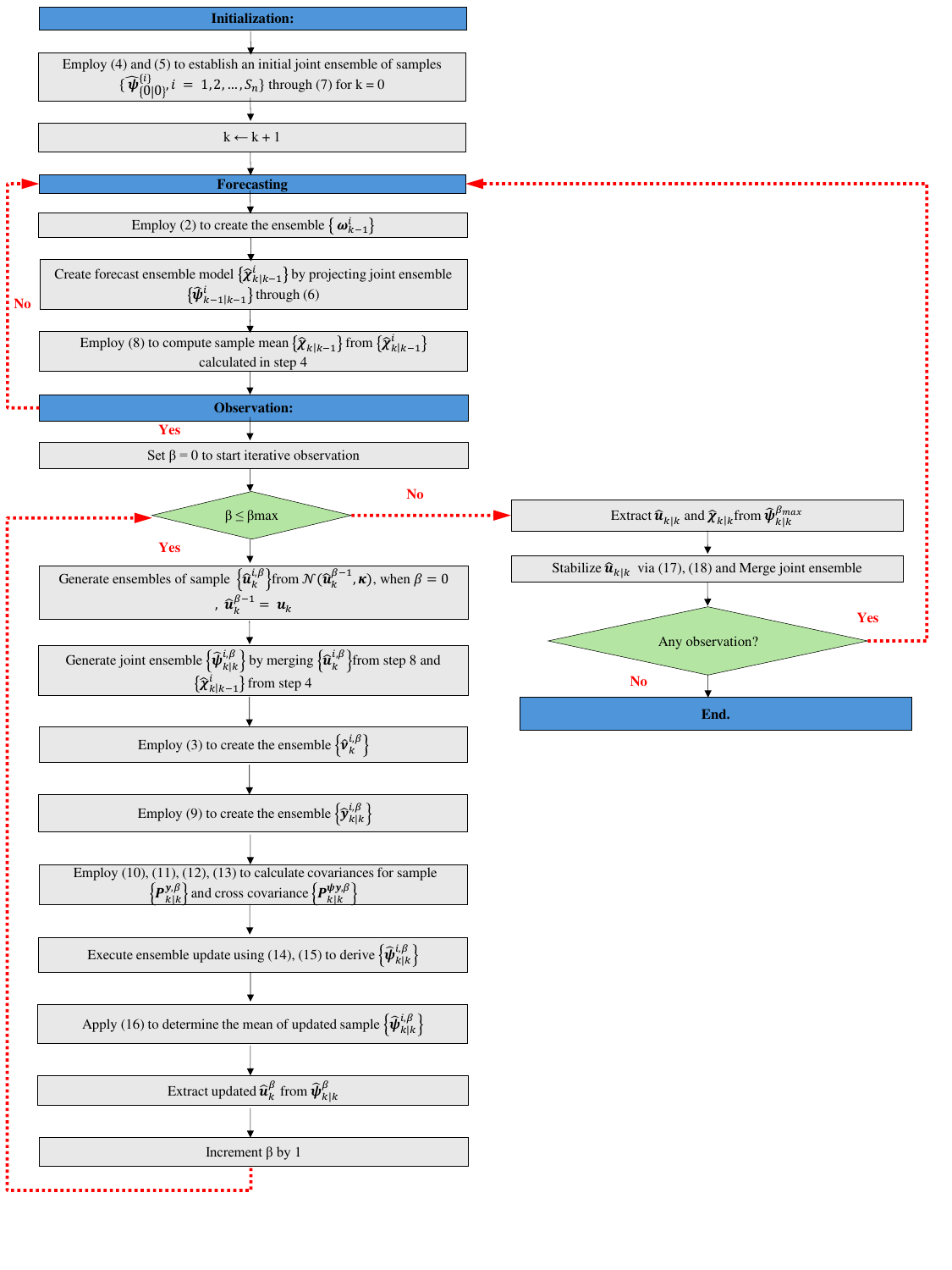} 
    \caption{The algorithm for EnSISF-wDF}
    \label{flowchart}
\end{figure}

The EnSISF-wDF algorithm, as illustrated in in Fig. \ref{flowchart} for the system outlined in Eq. \ref{eq:SISSFwdf} comprises two primary stages: predict and update. It utilizes ensembles that represents the conditional distributions of input parameters and system states in light of output measurements. Using this these ensembles, it estimates joint distributions based on sample means and covariances. To improve the estimation accuracy, the observation step is applied iteratively. This algorithm is particularly efficient and effective for high-dimensional linear and nonlinear systems, as it utilizes derivative-free computations, which enhance its applicability in complex scenarios \cite{Fang2017-tx}.

\subsection{Parameter Stabilization}\label{subsec:Parameter Constraint}
In the context of parameter estimation using EnKF, constraining parameters is a common practice to prevent divergence and ensure physically realistic estimates, as supported by the literature \cite{canuto2020ensemble, khaki2020calibrating}. Such constraints help maintain the stability and accuracy of the estimation process by keeping parameters within plausible bounds informed by physical data.

In our approach, we apply this conventional practice to the estimation of the parameter, i.e. inlet boundary condition within the abdominal aorta. 
Specifically, the inlet velocity parameter is constrained such that it does not exceed or fall below 80\% of the average velocity measured at predefined sensor locations immediately upstream of the inlet boundary. Let \( \mathbf{\hat{u}}_{inlet,n} \) denote the estimated parameter (inlet velocity condition) at the nth time step. The average velocity of the measurements from the sensors in front of the inlet is given by:
\begin{equation}
\label{eq_const1}
\bar{v}_m = \frac{1}{M} \sum_{j=1}^{M} v_{j,m}
\end{equation}

where \( v_{j,m} \) represents the velocity measured at the jth sensor at time step m, and M is the total number of sensor points. To enforce the constraint, the estimated parameter \( \mathbf{\hat{u}}_{inlet,n} \) is bounded by:
\begin{equation}
\label{eq_const2}
0.8 \cdot \bar{v}_m \leq \mathbf{\hat{u}}_{inlet,n} \leq 1.2 \cdot \bar{v}_m
\end{equation}

This approach ensures that the estimated inlet boundary condition stays within 80\% of the values measured at the stabilization points, indicated by purple circles in Figs. \ref{fig:2DInVivoFlowSensors} and \ref{fig:3DInVivoFlowSensors}.

\section{Synthetic Measurement Data}\label{sec:Dynamic Model}
To generate measurement data, one could use either in vivo measurements from 4D MRI or in vitro data. However, this research primarily aims to assess the effectiveness of the EnSISF for predicting boundary conditions of the aorta. Therefore, to produce synthetic measurement data—referred to here as 'pseudo-experimental data'—we employ high-fidelity numerical simulations that achieve a level of accuracy comparable to in vitro experimental data, providing a robust foundation for data assimilation. The specifics of this solver are detailed in this section.

\subsection{Mathematical Model}
The governing equations of the blood motion within a time-independent domain \( \Omega \) over a time interval \((t_0, t_F]\) are expressed as the momentum and mass conservation as follows: 
\begin{equation} 
\label{eq_NS} 
\begin{array}{ll} 
\rho \frac{\partial U_i}{\partial t} + \rho \frac{\partial}{\partial x_j} (U_i U_j) - \frac{\partial \sigma_{ij}}{\partial x_j} = 0 \quad \text{in } \Omega \times (t_0, t_F],
\end{array} 
\end{equation} 

\begin{equation} 
\label{eq_cont_cap}
\begin{array}{ll} 
\frac{\partial U_i}{\partial x_i} = 0 \quad \text{in } \Omega \times (t_0, t_F],
\end{array} 
\end{equation}
where $\rho=1060 \, \text{kg/m}^3$, $U_i$ and $\sigma_{ij}$ denote the blood density, velocity vector, and the Cauchy stress tensor. For both Newtonian and non-Newtonian fluids we have:
\begin{equation}
\label{eq_cauchy}
\begin{array}{ll}
\sigma_{ij}(U, p) = -p \delta_{ij} + 2 \mu(\dot{\gamma}) D_{ij},
\end{array}
\end{equation}

where $D_{ij} = \frac{1}{2} \left( \frac{\partial U_i}{\partial x_j} + \frac{\partial U_j}{\partial x_i} \right)$, \(p\) represents the pressure, $\mu$ is the dynamic viscosity and $\dot{\gamma}$ is the shear rate:
\begin{equation}
\label{eq_shearRate}
\begin{aligned}
\dot{\gamma} = \sqrt{\frac{1}{2} D_{ij} D_{ij}}.
\end{aligned}
\end{equation}

Under Newtonian blood assumption, the viscosity $\mu$ remains constant at 0.0035 Pa.s and does not depend on the shear rate. Conversely, the viscosity of a non-Newtonian fluid varies with the shear rate, $\dot{\gamma}$. The Casson model used in this simulation accurately represents the non-Newtonian characteristics of blood viscosity and is expressed as follows \citep{salavatidezfouli2024effect}:
\begin{equation}
\label{eq_Casson}
\mu=\frac{\boldsymbol{\tau}_0}{\dot{\gamma}} +\frac{\sqrt{\mu_{\infty} \boldsymbol{\tau}_0}}{\sqrt{\dot{\gamma}}}+\mu_{\infty},
\end{equation}
where $\boldsymbol{\tau}_0$ is 0.005 and $\mu_{\infty}$ is 0.0035.

\subsubsection{Turbulence}
The flow within the abdominal aorta (AA) is typically laminar during rest but becomes turbulent at peak systole \citep{etli2021numerical}. For a precise solution, neither purely laminar nor fully turbulent models are exact enough \citep{kelly2020influence}. However, the transitional $k-\omega$ SST model enhances the accuracy in such situations by incorporating equations for intermittency, $\gamma$, and transition momentum thickness Reynolds number, $\widetilde{Re}_{\theta t}$. Owing to its effectiveness, the transitional $k-\omega$ SST model is widely used in cardiovascular research to model blood flow in the aorta \citep{zhu2022fluid, sengupta2023aortic, salavatidezfouli2024effect}, making it the preferred choice for our simulations of generating pseudo-experimental data. It incorporates the following set of equations \citep{almohammadi2015modeling}:
\begin{equation}
\label{eq_kwSSTTransition}
\begin{aligned}
& \frac{\partial(\rho k)}{\partial t}+\frac{\partial\left(\rho U_i k \right)}{\partial x_i}=\widetilde{P_k}-E_k+\frac{\partial}{\partial x_j}\left[\Gamma_k \frac{\partial k}{\partial x_j}\right], \\
& \frac{\partial(\rho \omega)}{\partial t}+\frac{\partial\left(\rho U_i \omega\right)}{\partial x_i}=P_\omega-E_\omega+D_\omega+\frac{\partial}{\partial x_j}\left[\Gamma_\omega \frac{\partial \omega}{\partial x_j}\right], \\
& \frac{\partial(\rho \gamma)}{\partial t}+\frac{\partial(\rho U_i \gamma)}{\partial x_i}=P_{\gamma 1}-E_{\gamma 1}+P_{\gamma 2}-E_{\gamma 2}+\frac{\partial}{\partial x_j}\left[\left(\mu+\frac{\mu_t}{\sigma_\gamma}\right) \frac{\partial \gamma}{\partial x_j}\right], \\
& \frac{\partial(\rho \widetilde{Re}_{\theta t})}{\partial t}+\frac{\partial(\rho U_i \widetilde{Re}_{\theta t})}{\partial x_i}=P_{\theta t}+\frac{\partial}{\partial x_j}\left[\left(\mu+\mu_t\right) \frac{\partial \widetilde{Re}_{\theta t}}{\partial x_j}\right].
\end{aligned}
\end{equation}

The production terms for kinetic energy, specific dissipation, and transition sources are represented by \(\tilde{P_k}\), \(P_\omega\), \(P_{\gamma_1}\), and \(E_{\gamma_1}\), respectively. Similarly, the destruction terms of kinetic energy, specific dissipation, and transition sources are denoted by \(E_k\), \(E_\omega\), \(E_{\gamma_2}\), and \(E_{\gamma_2}\), respectively. \(D_\omega\) represents the cross-diffusion term, and \(P_{\theta t}\) is the source term for $\widetilde{Re}_{\theta t}$ \cite{johari2019disturbed}.

\subsubsection{Boundary Condition}\label{section:BC}
\begin{figure}[H]
  \centering
  \includegraphics[width=0.75\textwidth]{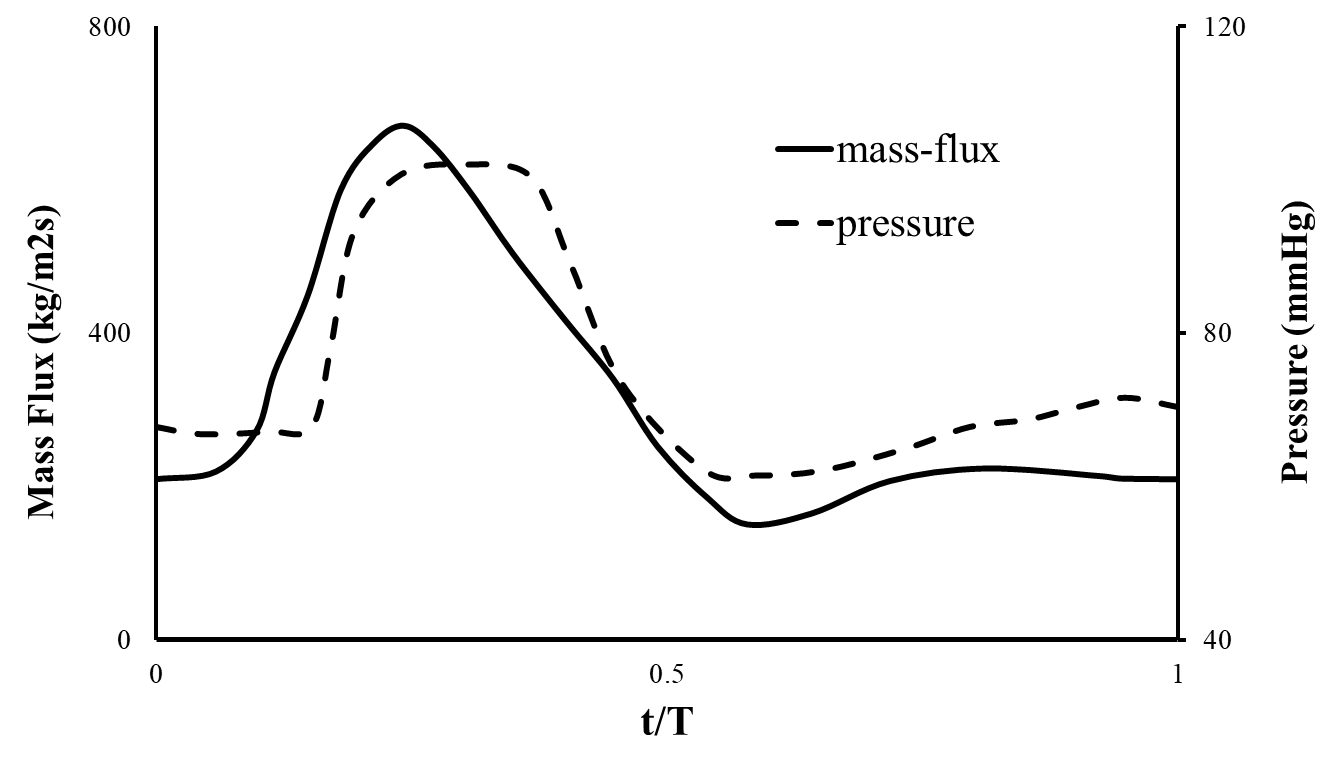}
  \caption{Temporal profiles for the inlet and outlet boundary conditions.}
  \label{FIG_BC}
\end{figure}

To generate pseudo-experimental data, accurately defined boundary conditions are essential at both the inlet and outlet. The inflow and outflow boundaries are set using mass flux inlet and pressure outlet conditions, respectively, based on data from \cite{alishahi2011numerical}. The time-dependent profiles for these inlet and outlet conditions are shown in Fig. \ref{FIG_BC}.

\subsubsection{Numerical Solution}
ANSYS FLUENT, a commercial package based on the cell-centered finite volume approach, has been employed to model the blood flow in the aortic model. Both the spatial and temporal terms in Eq. \ref{eq_NS} are discretized using a second-order scheme, and the SIMPLE (Semi-Implicit Method for Pressure-Linked Equations) algorithm has been incorporated for pressure-velocity coupling.

\subsection{2D Ideal Model}
\begin{figure}[H]
    \centering
    \begin{subfigure}[b]{0.475\textwidth}
        \centering
        \includegraphics[width=\textwidth]{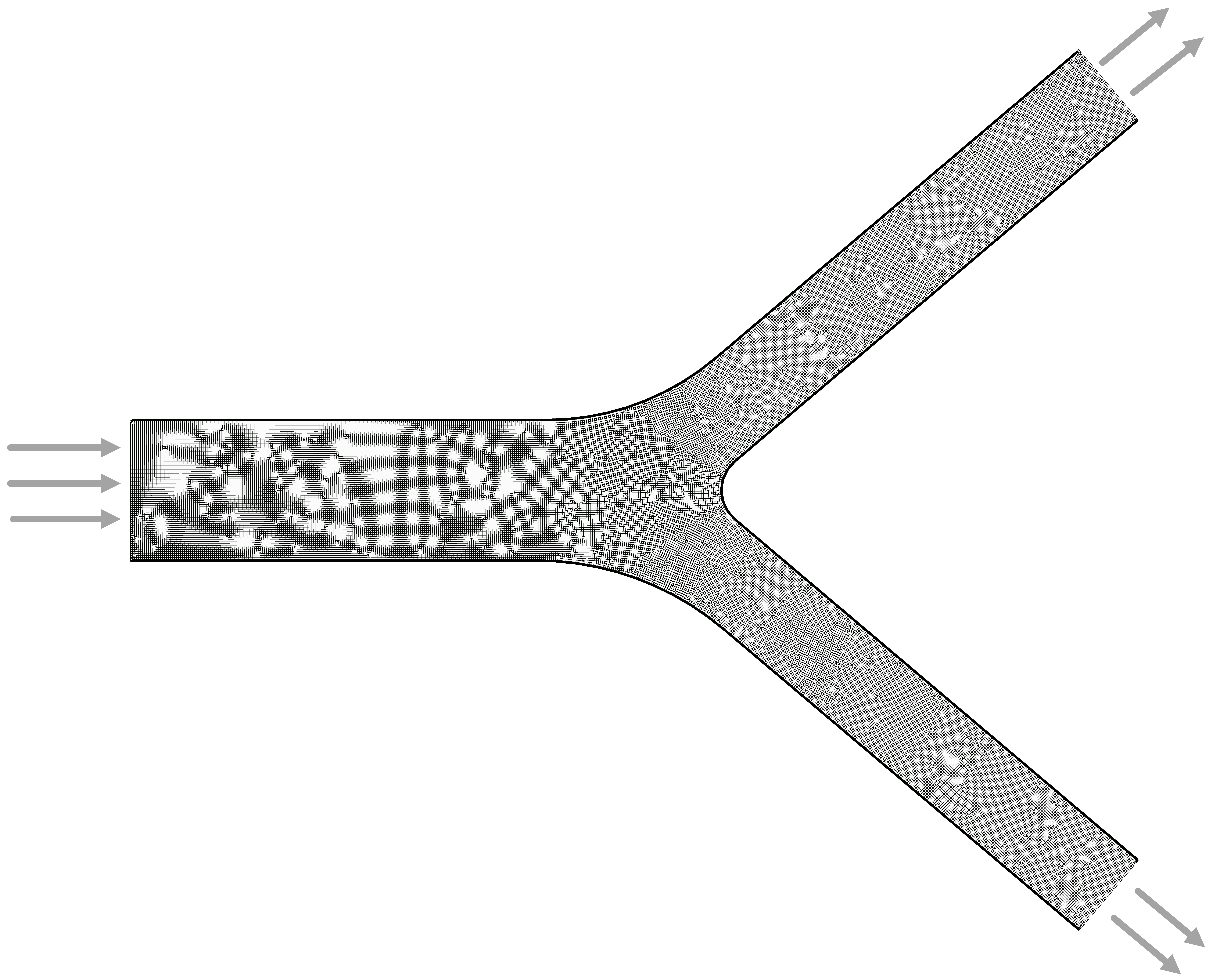}
        \caption{Fine mesh}
        \label{fig:mesh-a}
    \end{subfigure}
    \hfill
    \begin{subfigure}[b]{0.475\textwidth}
        \centering
        \includegraphics[width=\textwidth]{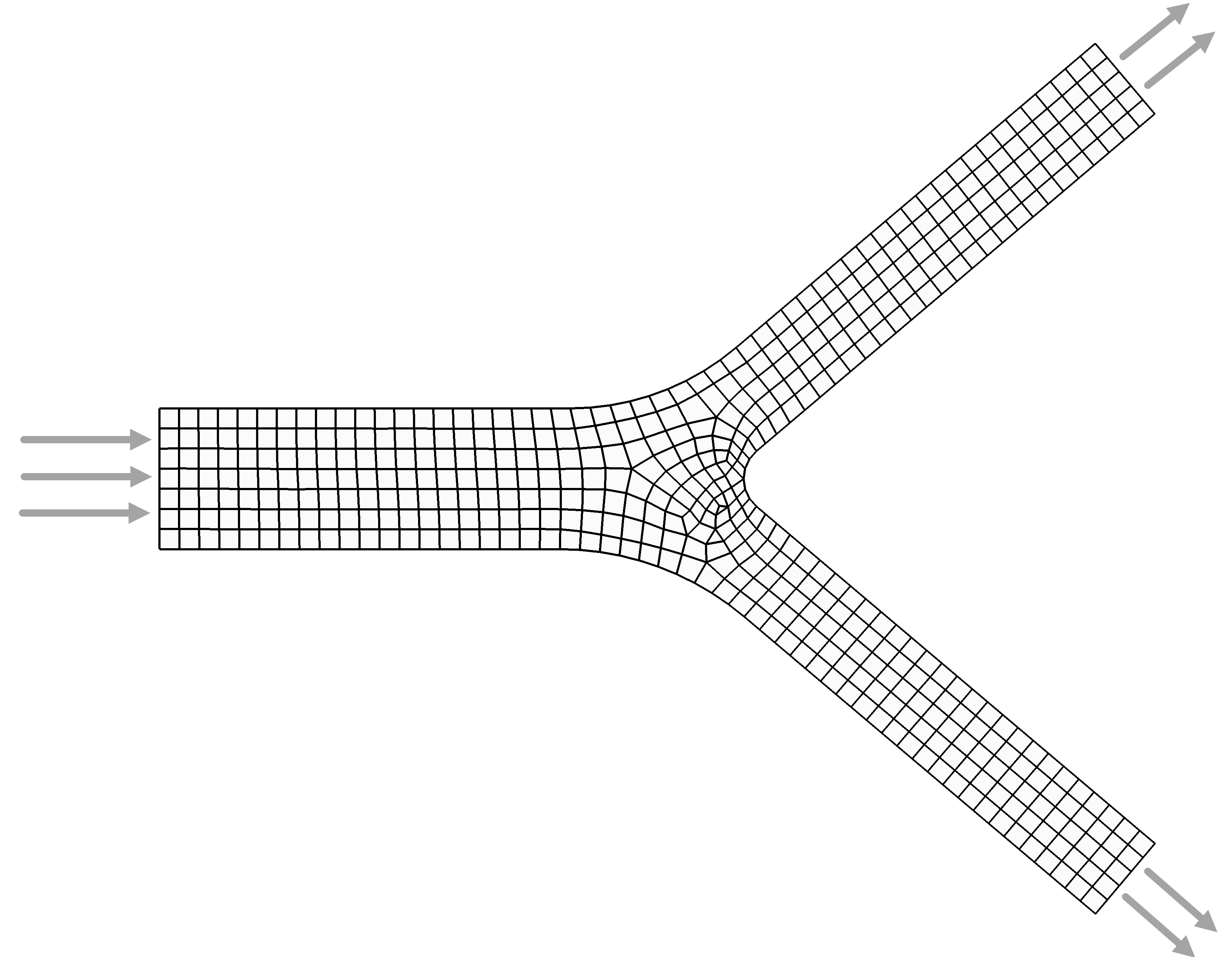}
        \caption{Coarse Mesh}
        \label{fig:mesh-b}
    \end{subfigure}
    
    \caption{Demonstration of 2D computational meshes for a)generation of pseudo-experimental data and b)forward solver problem.}

    \label{FIG_meshes2d}
\end{figure}

To obtain high-fidelity pseudo-experimental data, a fine mesh with 8000 elements—primarily quadrilateral—was created based on a mesh independence study. This mesh includes five boundary layers near the wall boundaries and features a very fine distribution throughout the domain, enabling precise flow and gradient resolution within the aorta. The mesh is illustrated in Fig.  \ref{fig:mesh-a}.

\subsection{3D Patient-Specific Model}
\begin{figure}[H]
    \centering
    \begin{subfigure}[b]{0.475\textwidth}
        \centering
        \includegraphics[width=\textwidth]{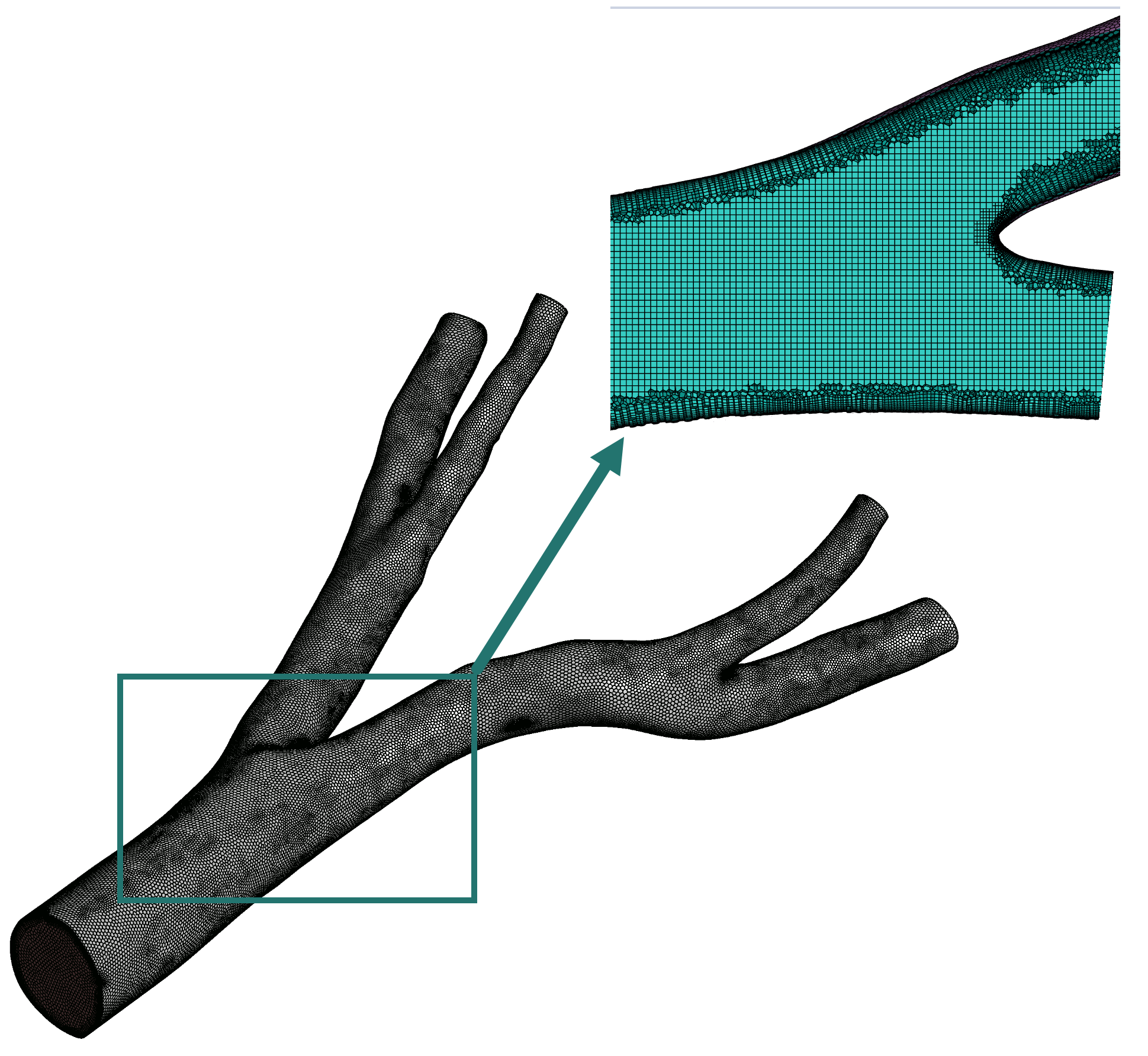}
        \caption{Fine mesh}
        \label{fig:mesh3d-a}
    \end{subfigure}
    \hfill
    \begin{subfigure}[b]{0.475\textwidth}
        \centering
        \includegraphics[width=\textwidth]{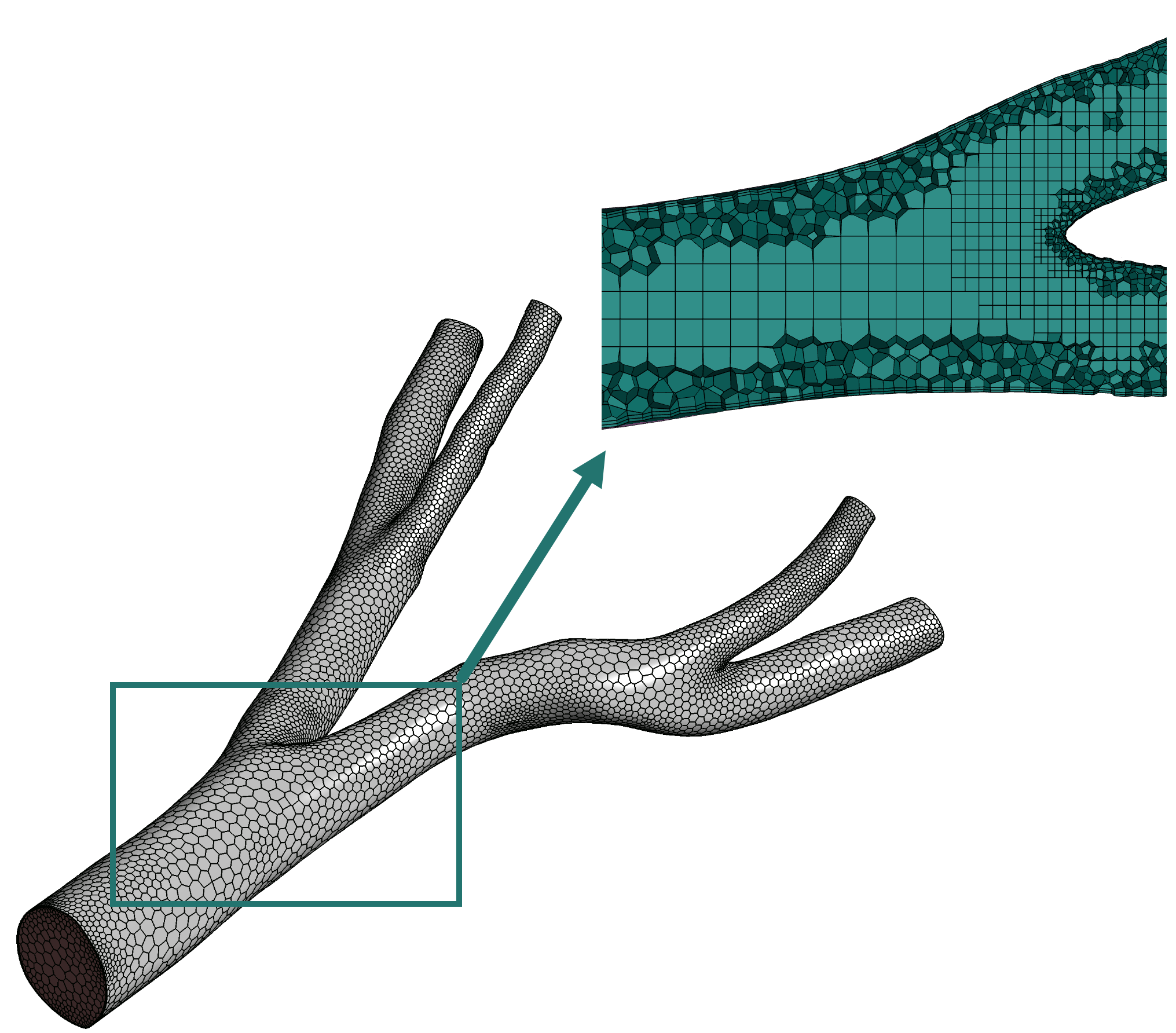}
        \caption{Coarse Mesh}
        \label{fig:mesh3d-b}
    \end{subfigure}
    
    \caption{Demonstration of 3D meshes for a)generation of pseudo-experimental data and b)forward solver problem.}

    \label{FIG_meshes3d}
\end{figure}

Finally, to extend the generality of our ensemble model, we applied it to a 3D model specific to a patient’s abdominal aorta. Similar to the 2D case, a fine mesh is designed to produce highly accurate pseudo-experimental data consisting of $1.5\times10^5$ elements including 10 boundary layers, ensuring an acceptable y-plus range ($y^+ \sim 1$). The mesh, along with a cross-section, is illustrated in Fig. \ref{fig:mesh3d-a}.

\section{Forward Solver}
It is important to note that the efficacy of the EnKF employed in this study does not rely on the precision of the underlying physical model to perform corrections on the flow field or boundary condition predictions. Instead, the EnSISF's strength lies in its ability to integrate and adjust based on incoming data, irrespective of the model's fidelity or any prior knowledge of the boundary conditions. 

Thus, the model employed during the forecast loop (as outlined in the flowchart shown in Fig. \ref{flowchart}) utilizes a lower-fidelity approach. This model simulates blood flow using a laminar flow assumption and Newtonian viscosity, computed over coarser meshes for 2D and 3D cases, as illustrated in Figs. \ref{fig:mesh-b} and \ref{fig:mesh3d-b}. This choice of a less complex model is strategic, enabling faster computations as the forecast model runs multiple simulations across numerous ensembles, denoted by $S_n$. Using a high-fidelity model in this context would be prohibitively time-consuming and unnecessary for the iterative updates performed by this advanced extension of ensemble Kalman filter.

\subsection{Measurement Points}
The EnKF model requires several measurement points within the domain to improve accuracy in state and parameter prediction. The number of these measurement points influences prediction accuracy, with \cite{bakhshaei2024optimized} suggesting that approximately 5\% of total data points be used as measurement data. Accordingly, this study selected 27 out of 498 cells for the 2D case and 330 out of 8471 cells for the 3D case as measurement points, where high-fidelity pseudo-experimental data will be assimilated.

The measurement points—referred to as Flow Sensors—are strategically distributed across various sections of the aorta, particularly near bifurcation areas where flow dynamics are complex.

\begin{figure}[H]
    \centering
    \includegraphics[width=0.8\textwidth]{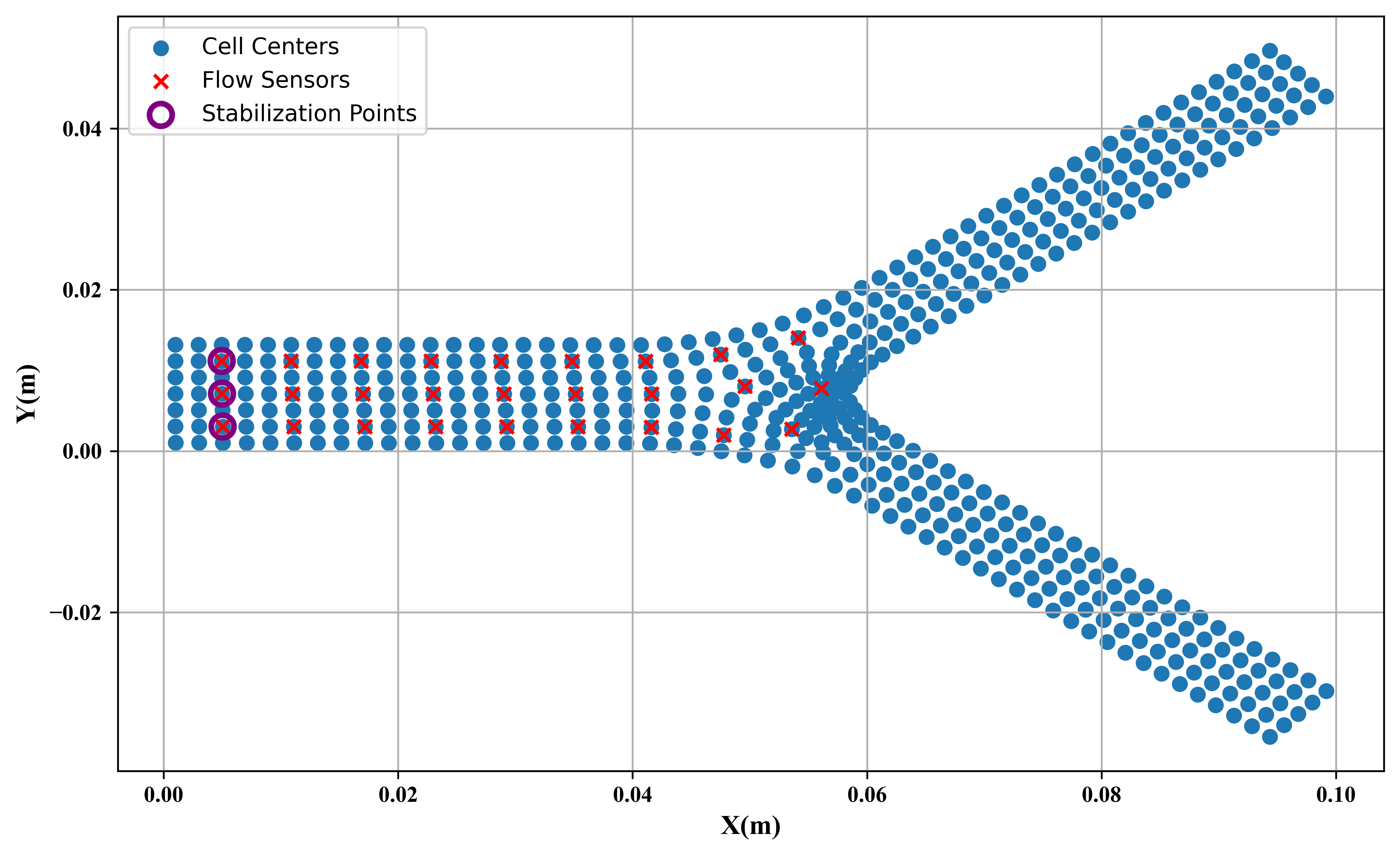}
    \caption{Spatial distribution of cell centers and flow sensors, stabilization points within the vascular model.}
    \label{fig:2DInVivoFlowSensors}
\end{figure}

\begin{figure}[H]
    \centering
    \includegraphics[width=1\linewidth]{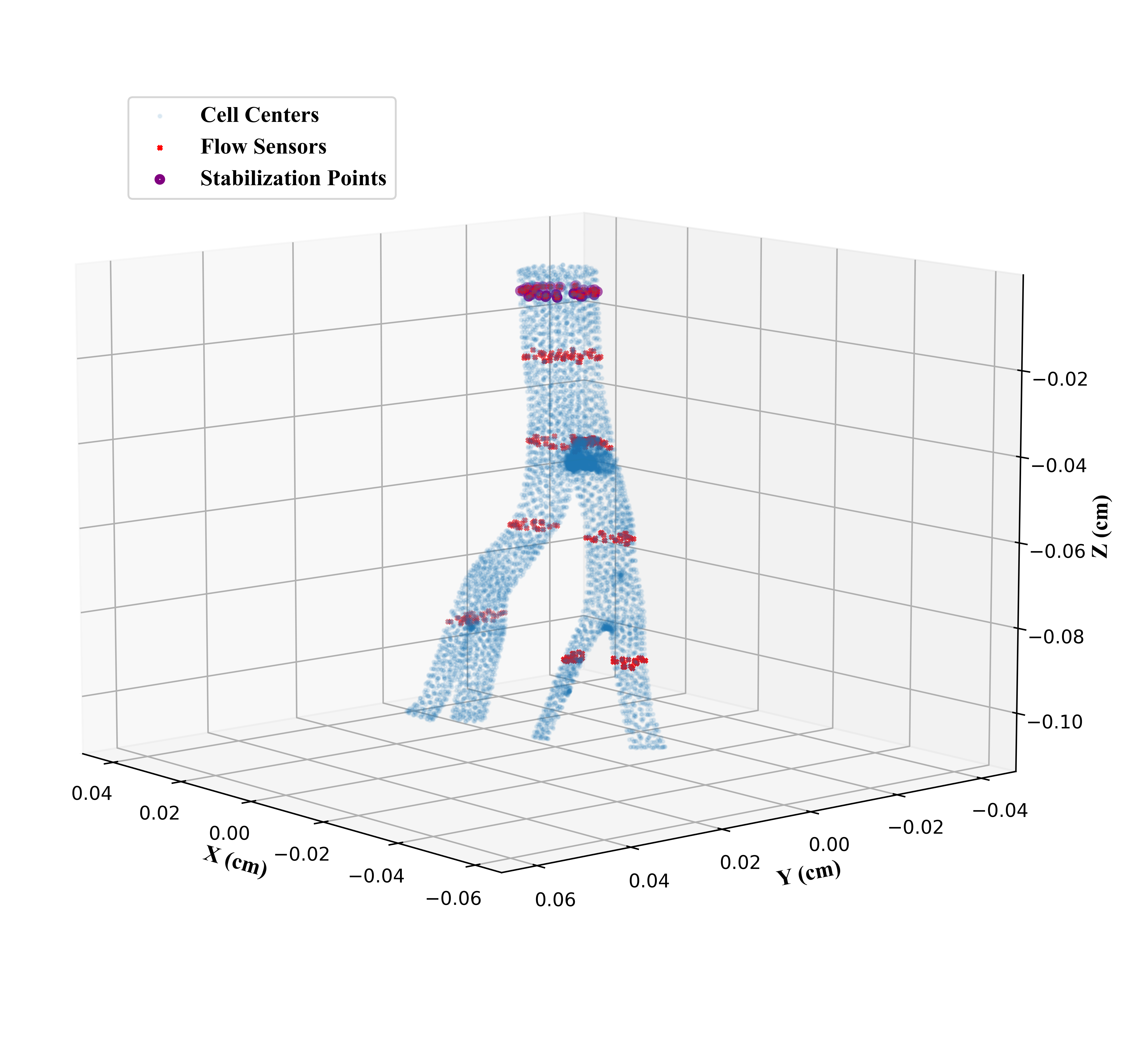}
    \caption{3D distribution of cell centers and flow sensors, and stabilization points within the patient-specific abdominal aorta model.}
    \label{fig:3DInVivoFlowSensors}
\end{figure}

Figs. \ref{fig:2DInVivoFlowSensors} and \ref{fig:3DInVivoFlowSensors} show the spatial arrangement of computational cell centers (blue dots) and flow sensor placements (red crosses) for both the 2D ideal model and the 3D patient-specific model. These cell centers align with nodes in the coarse mesh used during the ensemble forecast loop.

%

\section{Results and Discussion}\label{sec:resVessel}
This section discusses parameter prediction for both 2D ideal and 3D patient-specific models. To develop an optimized solution strategy, we evaluated our EnKF model with increasing levels of complexity. Starting with steady parameter prediction for the 2D model, we assume a constant boundary condition over time, with results detailed in Section \ref{RD:2DConstant}. Next, we introduce time-dependency by using a spatially uniform but temporally variable boundary condition for the 2D model, as outlined in Section \ref{RD:2DTime}. Finally, we consider a parameter that varies with both time and space for the 2D model, with results in Section \ref{RD:2DTimeSpace}. After analyzing these progressively complex cases for the 2D ideal model, we apply the EnKF to the 3D patient-specific aorta model with a time-space-dependent boundary condition, as shown in Section \ref{3DTimeSpace}.

\subsection{2D Ideal Model}\label{RD:2DModels}
For the 2D case, after using a trial-and-error approach to optimize the model hyperparameters, we selected the parameter set shown in Table \ref{tab:optim_par_2D}, which achieved the highest accuracy in parameter estimation across different scenarios.

\begin{table}[H]
\caption{Optimal parameters for three different scenarios of 2D ideal model}
\label{tab:optim_par_2D}
\resizebox{\textwidth}{!}{
\begin{tabular}{|
>{\columncolor[HTML]{C0C0C0}}l |
>{\columncolor[HTML]{FFFFFF}}c |
>{\columncolor[HTML]{FFFFFF}}c |
>{\columncolor[HTML]{FFFFFF}}c |}
\hline
\cellcolor[HTML]{FFFFFF}\textbf{Variables/Scenario} & \cellcolor[HTML]{C0C0C0}\textbf{\begin{tabular}[c]{@{}c@{}}Constant \\ Parameter\end{tabular}} & \cellcolor[HTML]{C0C0C0}\textbf{\begin{tabular}[c]{@{}c@{}}Time Dependent\\  Parameter\end{tabular}} & \cellcolor[HTML]{C0C0C0}\textbf{\begin{tabular}[c]{@{}c@{}}Space-time Dependent\\  Parameter\end{tabular}} \\ \hline
\textbf{Number of Seeds}                            & {\color[HTML]{343434} 80}                                                                      & 80                                                                                                   & 80                                                                                                         \\ \hline
\textbf{Prior State Mean}                           & 0.01                                                                                           & 0.1                                                                                                  & 0.1                                                                                                        \\ \hline
\textbf{Prior State Covariance}                     & 1E-10                                                                                          & 1E-04                                                                                                & 1E-04                                                                                                      \\ \hline
\textbf{Model Error Covariance (Qk)}                & 1E-08                                                                                          & 1E-04                                                                                                & 1E-04                                                                                                      \\ \hline
\textbf{Measurements Noise Covariance (Rk)}         & 1E-10                                                                                          & 1E-08                                                                                                & 1E-08                                                                                                      \\ \hline
\textbf{Prior Parameter Mean}                       & 0.015                                                                                          & 0.1                                                                                                  & 0.1                                                                                                        \\ \hline
\textbf{Prior Parameter Covariance}                 & 4E-06                                                                                          & 4E-04                                                                                                & 4E-04                                                                                                      \\ \hline
\textbf{Delta t (sec)}                              & 0.01                                                                                           & 0.01                                                                                                 & 0.01                                                                                                       \\ \hline
\textbf{Observation Span (sec)}                     & {[}0.02, 0.04, 0.05{]}                                                                         & 0.02                                                                                                 & 0.02                                                                                                       \\ \hline
\textbf{T final (sec)}                              & 1                                                                                              & 1                                                                                                    & 1                                                                                                          \\ \hline
\end{tabular}
}
\end{table}

\subsubsection{Constant Parameter}\label{RD:2DConstant}

The observation span significantly impacts the accuracy of the EnKF method. Table \ref{tab:error_constantParameter} shows the relative error associated with different observation spans. After each cycle of data assimilation, the average relative error of the parameter estimation is calculated as follows: 
\begin{center}
\begin{equation}
\text{Mean Relative Error} = \frac{1}{N} \sum_{i=1}^{N} \frac{\left| P_{\text{actual}, i} - P_{\text{predicted}, i} \right|}{\left| P_{\text{actual}, i} \right|}
\end{equation}
\end{center}

where \( N \) represents the total count of time steps through the cycle, i.e. 100.

\begin{figure}[H]
    \centering
    
    \begin{subfigure}{0.7\textwidth}
        \includegraphics[width=\linewidth]{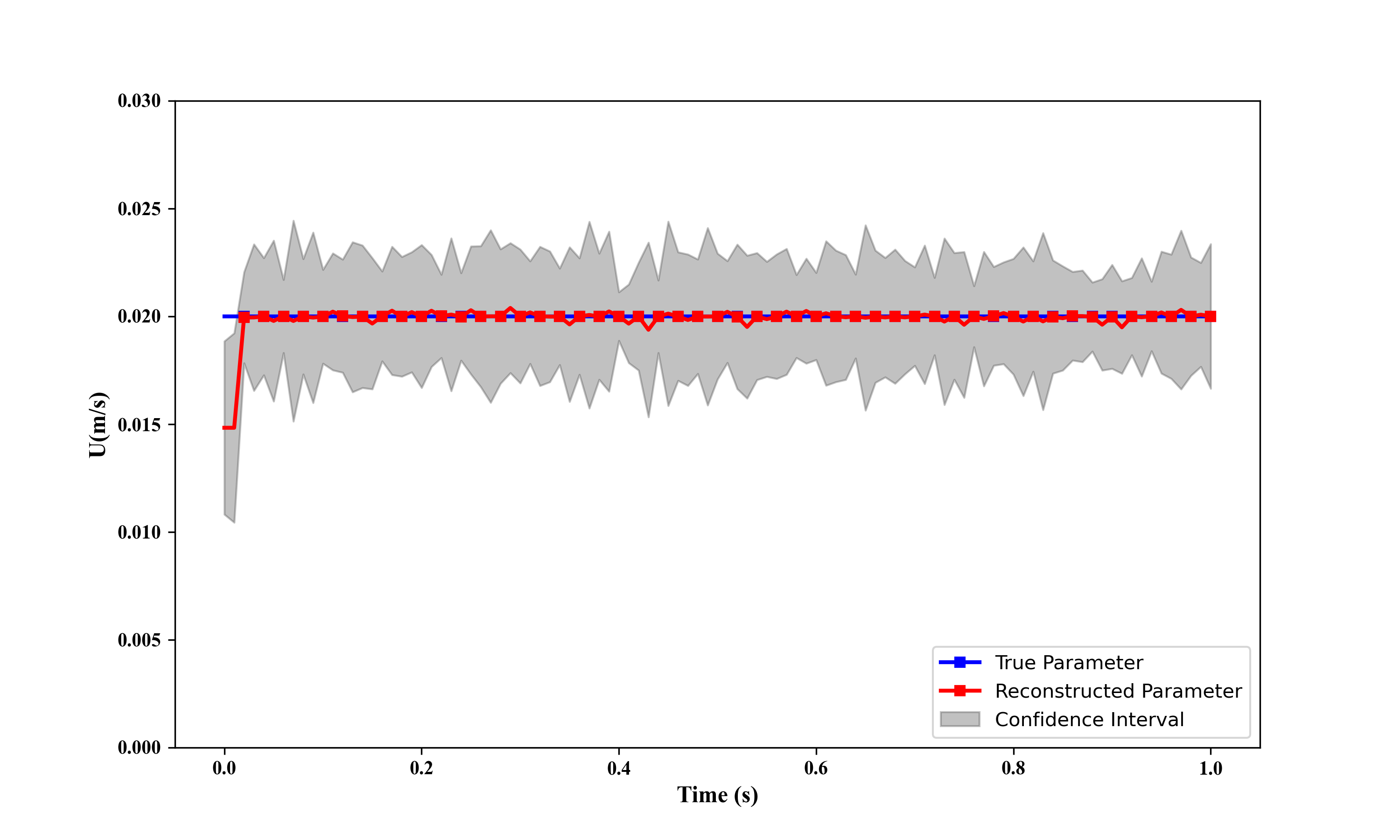}
        \caption{}
        \label{fig:ConstantParameter0.02}
    \end{subfigure}
    
    \begin{subfigure}{0.7\textwidth}
        \includegraphics[width=\linewidth]{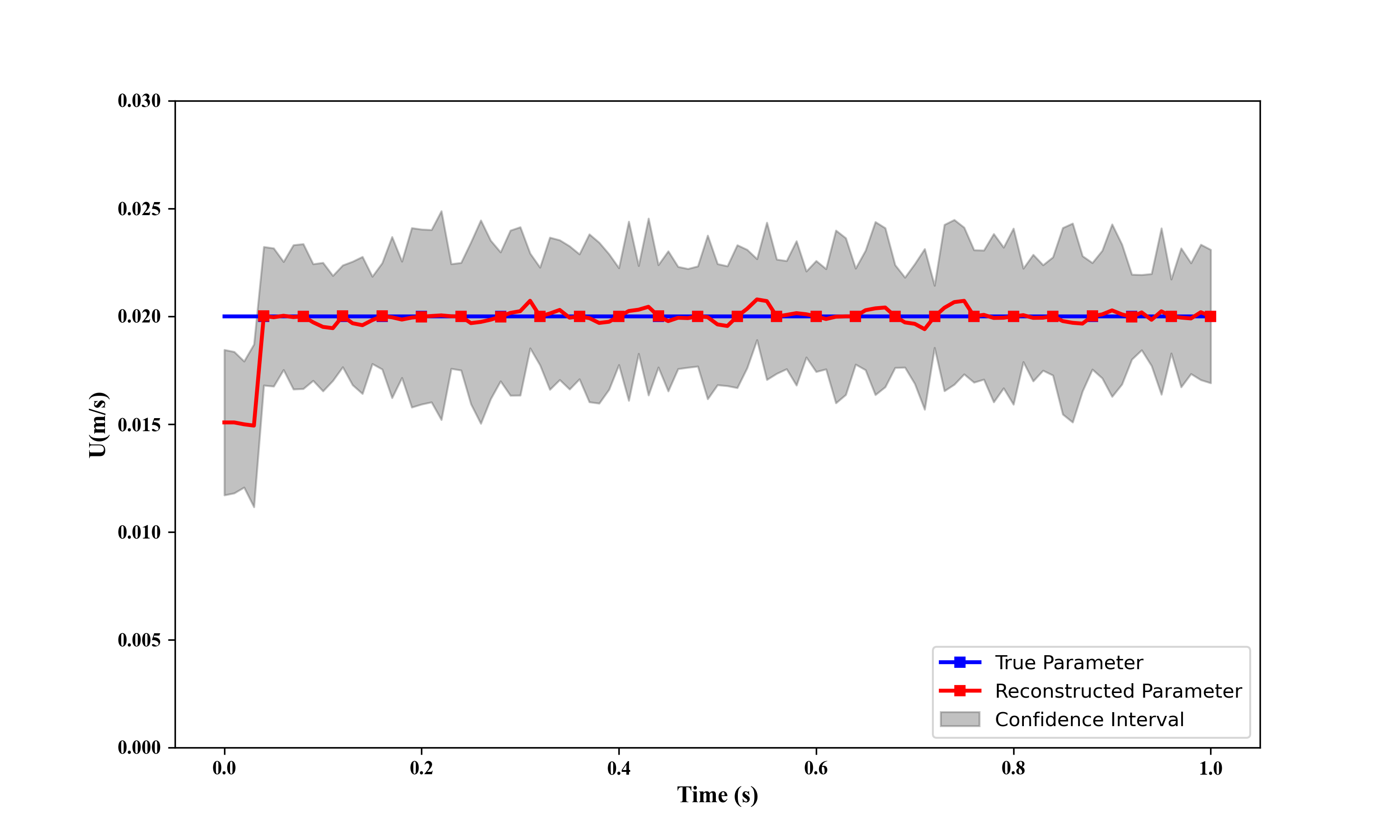} 
        \caption{}
        \label{fig:ConstantParameter0.04}
    \end{subfigure}
    
    \begin{subfigure}{0.7\textwidth}
        \includegraphics[width=\linewidth]{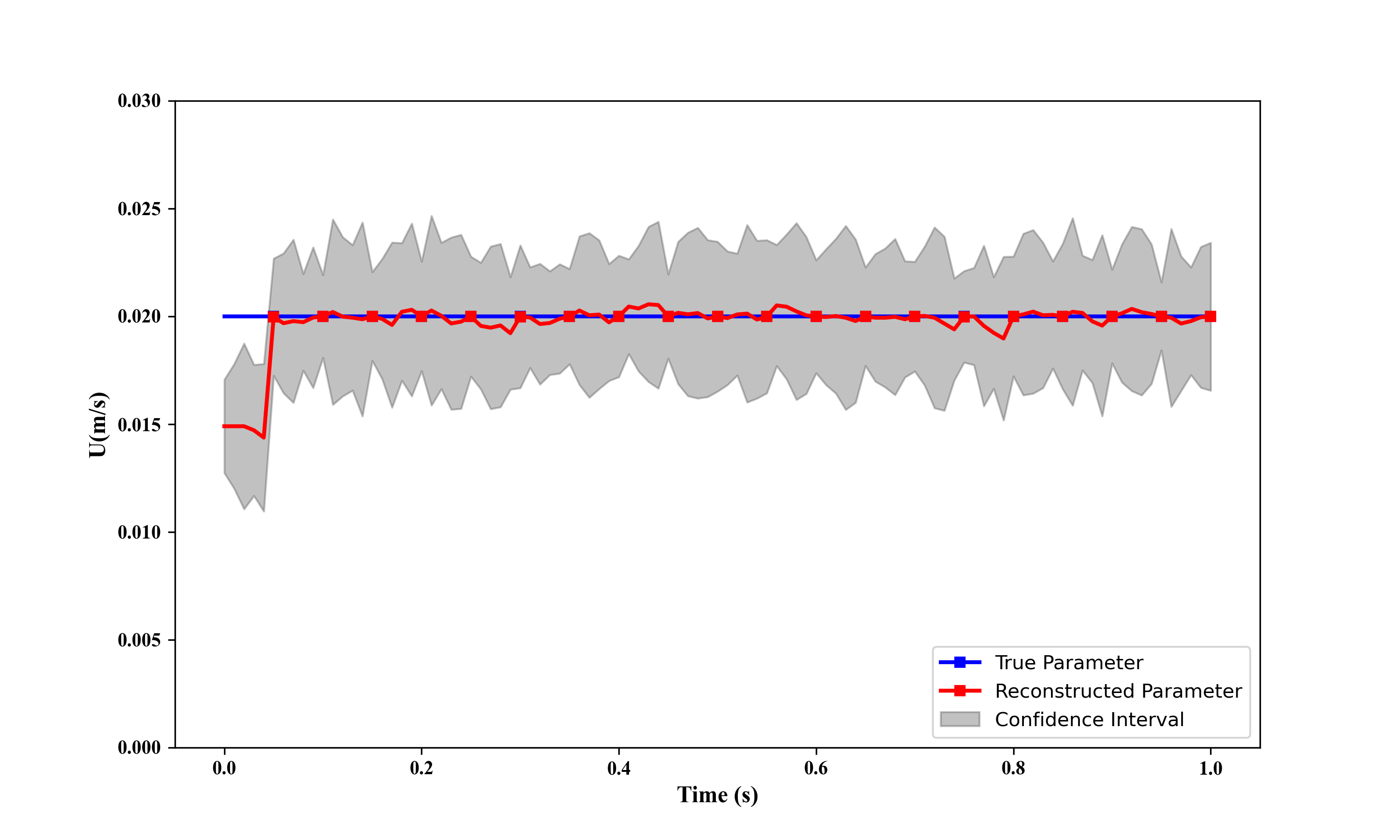} 
        \caption{}
        \label{fig:ConstantParameter0.05}
    \end{subfigure}
    
    \caption{Reconstruction of a constant parameter using varying observation spans. (a), (b), and (c) display the reconstruction results for observation spans of 0.02 sec, 0.04 sec, and 0.05 sec, respectively.}
    \label{fig:ConstantParameter}
\end{figure}

As indicated in Table \ref{tab:error_constantParameter}, increasing the observation span results in a higher error percentage. Accordingly, Fig. \ref{fig:ConstantParameter} illustrates the parameter reconstruction using EnSISF under those observation spans. The blue line represents the true parameter, the reconstructed parameter by the red line, and the grey-shaded area indicates the confidence interval. The square red markers indicate the times at which observations were made. Although the method shows a quick convergence of the reconstructed parameter to the true value within all subplots, an observation span of every two time steps has been chosen to achieve optimal accuracy for the rest of the simulations.
\begin{table}[H]
\centering
\caption{Relative error for constant parameter for three different observation spans}
\label{tab:error_constantParameter}
\begin{tabular}{|c|c|c|}
\hline
\rowcolor[HTML]{C0C0C0} 
\textbf{Observation Span (\#Time Step)} & \textbf{Observation Span (sec)} & \multicolumn{1}{l|}{\cellcolor[HTML]{C0C0C0}\textbf{Relative Error (\%)}} \\ \hline
2                  & 0.02                            & 0.996                                                                     \\ \hline
4                  & 0.04                            & 1.830                                                                     \\ \hline
5                  & 0.05                            & 2.180                                                                     \\ \hline
\end{tabular}
\end{table}
%

    
    
    
    
%

%

\begin{figure}[H]
    \centering
    
    \begin{subfigure}{0.7\textwidth}
        \includegraphics[width=\linewidth]{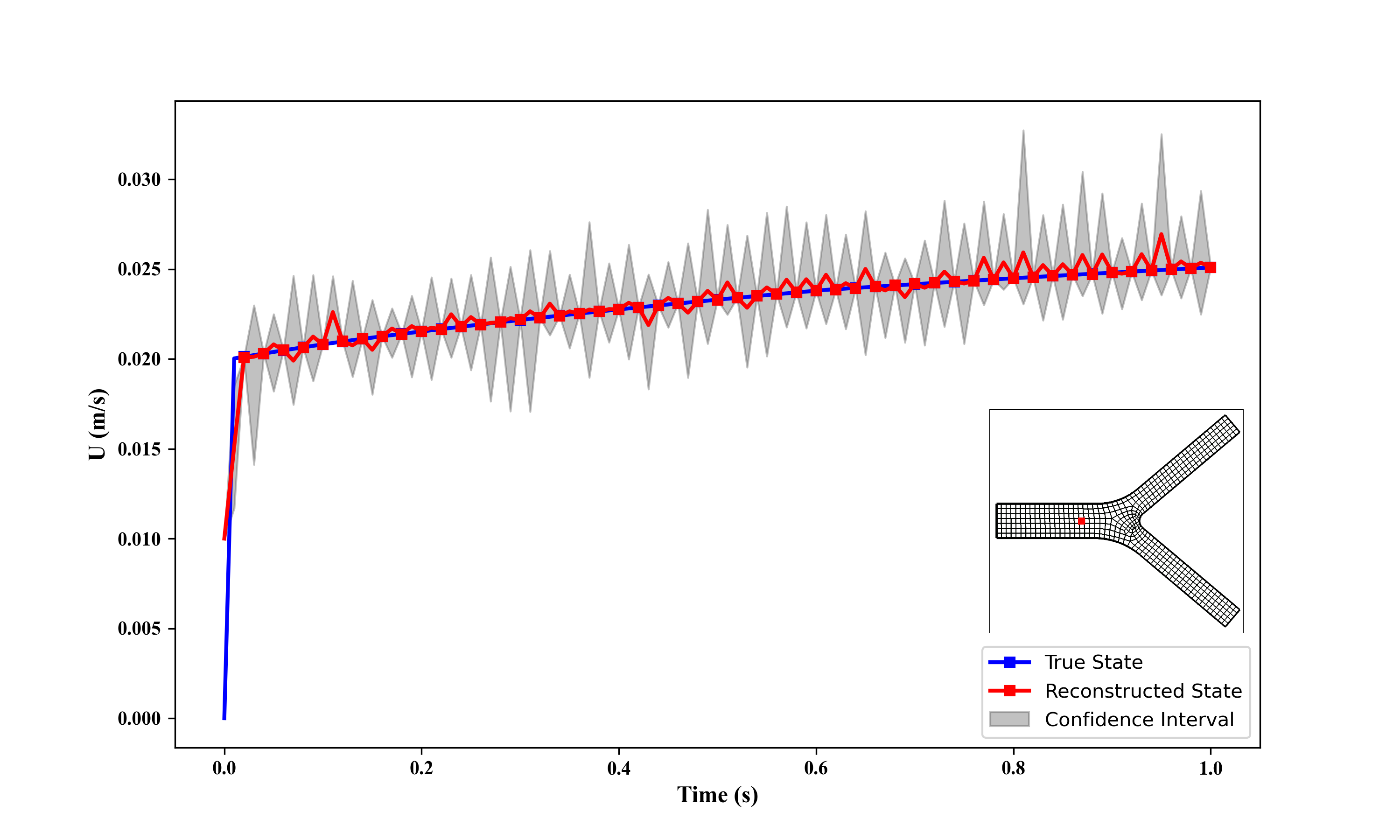}
        \caption{}
        \label{fig:State0.02}
    \end{subfigure}
    
    \begin{subfigure}{0.7\textwidth}
        \includegraphics[width=\linewidth]{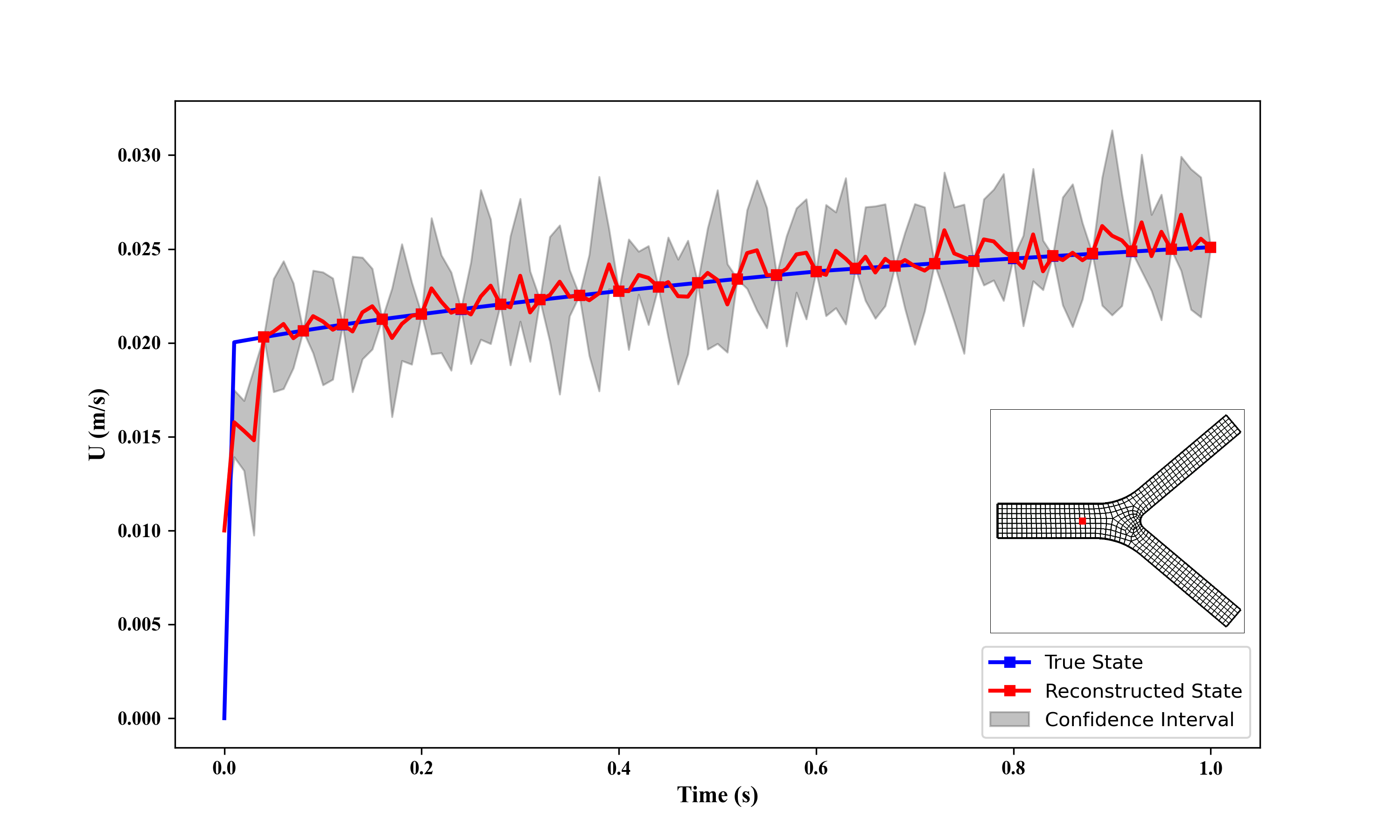} 
        \caption{}
        \label{fig:State0.04}
    \end{subfigure}
    
    \begin{subfigure}{0.7\textwidth}
        \includegraphics[width=\linewidth]{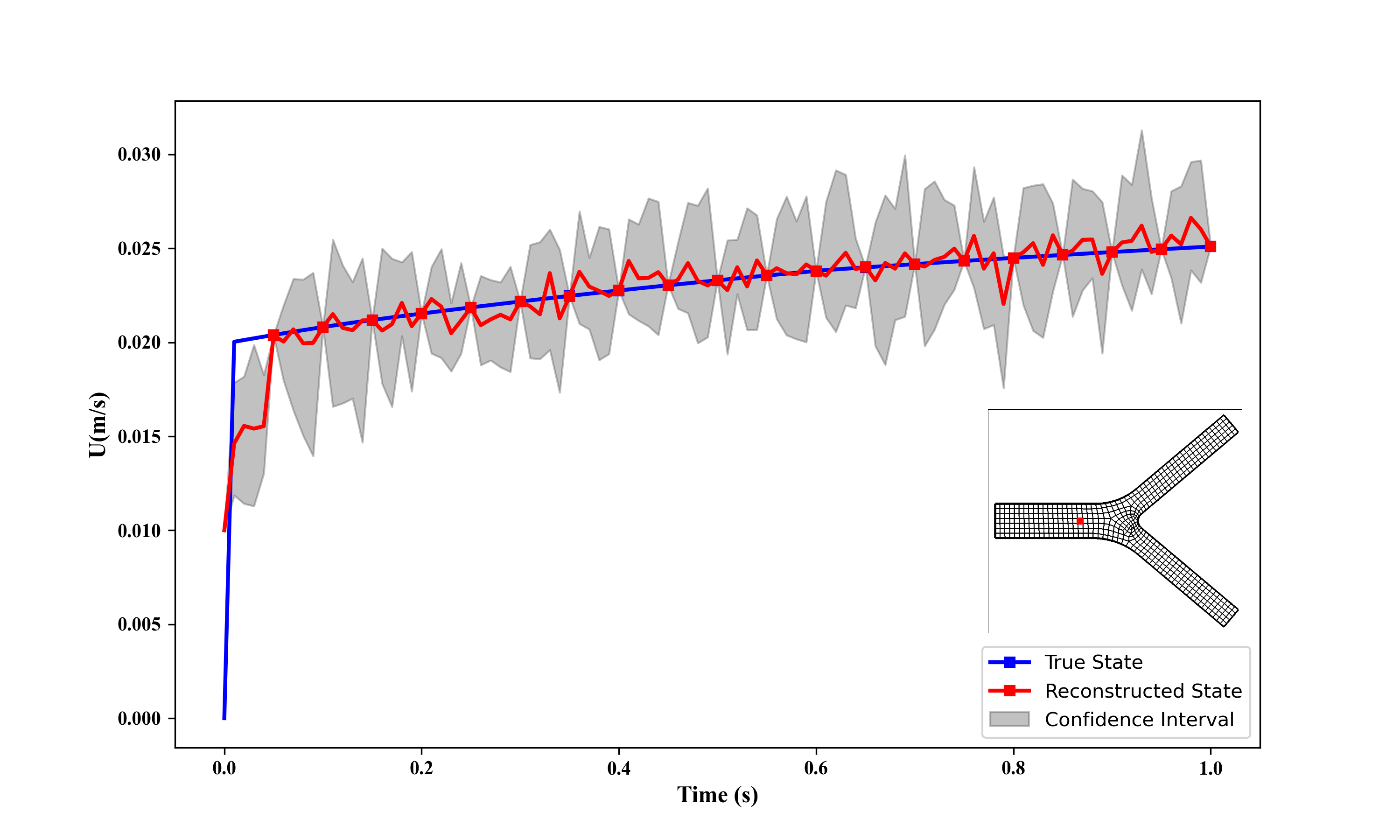} 
        \caption{}
        \label{fig:State0.05}
    \end{subfigure}
    
    \caption{Reconstruction of x-velocity state at varying observation intervals. Subplots (a), (b), and (c) show the reconstructed state for observation intervals of 0.02 sec, 0.04 sec, and 0.05 sec, respectively, at the measurement location ($x = 0.035147, y = 0.007069$).}
    \label{fig:stateForConstantParameter}
\end{figure}

Fig. \ref{fig:stateForConstantParameter} illustrates the reconstructed x-velocity state over time across varying observation intervals. The data is recorded at the measurement point located at $x = 0.035147$ and $y = 0.007069$, as indicated in each subplot. Initially, each subplot reveals a notable discrepancy between the true and reconstructed states, which narrows at the first observation interval, where data assimilation aligns the reconstructed state more closely with the true state. The confidence interval tightens significantly at observation points and widens between them, reflecting a dynamic adjustment of uncertainty as new data is incorporated. Extending the observation interval increases the confidence interval for times when data assimilation is not applied, i.e., between observation intervals. Nevertheless, the EnSISF method remains stable, showing robust performance without divergence during observation times.

%
%
%
%
\subsubsection{Time-Dependent Parameter}\label{RD:2DTime}

\begin{figure}[H] 
    \centering
\includegraphics[width=1\textwidth]{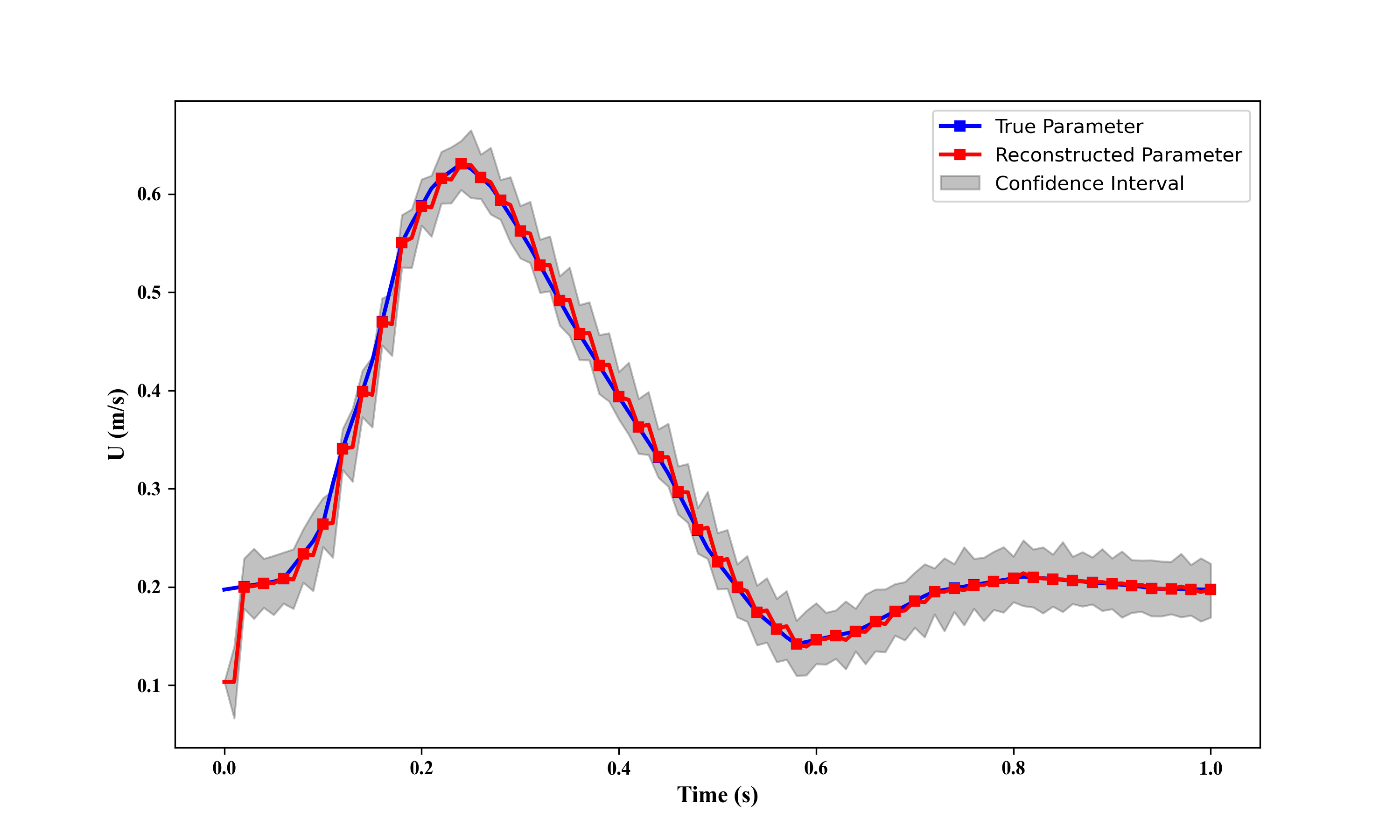}

    \caption{Comparison of true vs. reconstructed time-dependent velocity parameter}
    \label{fig:TimeDependentParameter}
\end{figure}

We now assess the EnSISF under the added complexity of a time-varying parameter, specifically a uniform inlet velocity that changes over time. Fig. \ref{fig:TimeDependentParameter} shows the reconstruction of this time-dependent velocity profile over a single cardiac cycle, with observations taken every 0.02 seconds. In the first two time steps, a noticeable discrepancy appears between the true and reconstructed parameters due to the lack of measurement assimilation. As observations are incorporated, the reconstructed parameter progressively aligns with the true parameter. The confidence interval narrows at observation points, reflecting reduced uncertainty with new data integration, and broadens between observations, indicating increased uncertainty. However, compared to the simpler case of a constant parameter over time, the confidence interval of the reconstructed parameter does not narrow as much, even at observation points, due to the non-linearity and complexity inherent in the cardiovascular model.

\begin{figure}[H]
    \centering
\includegraphics[width=1\textwidth]{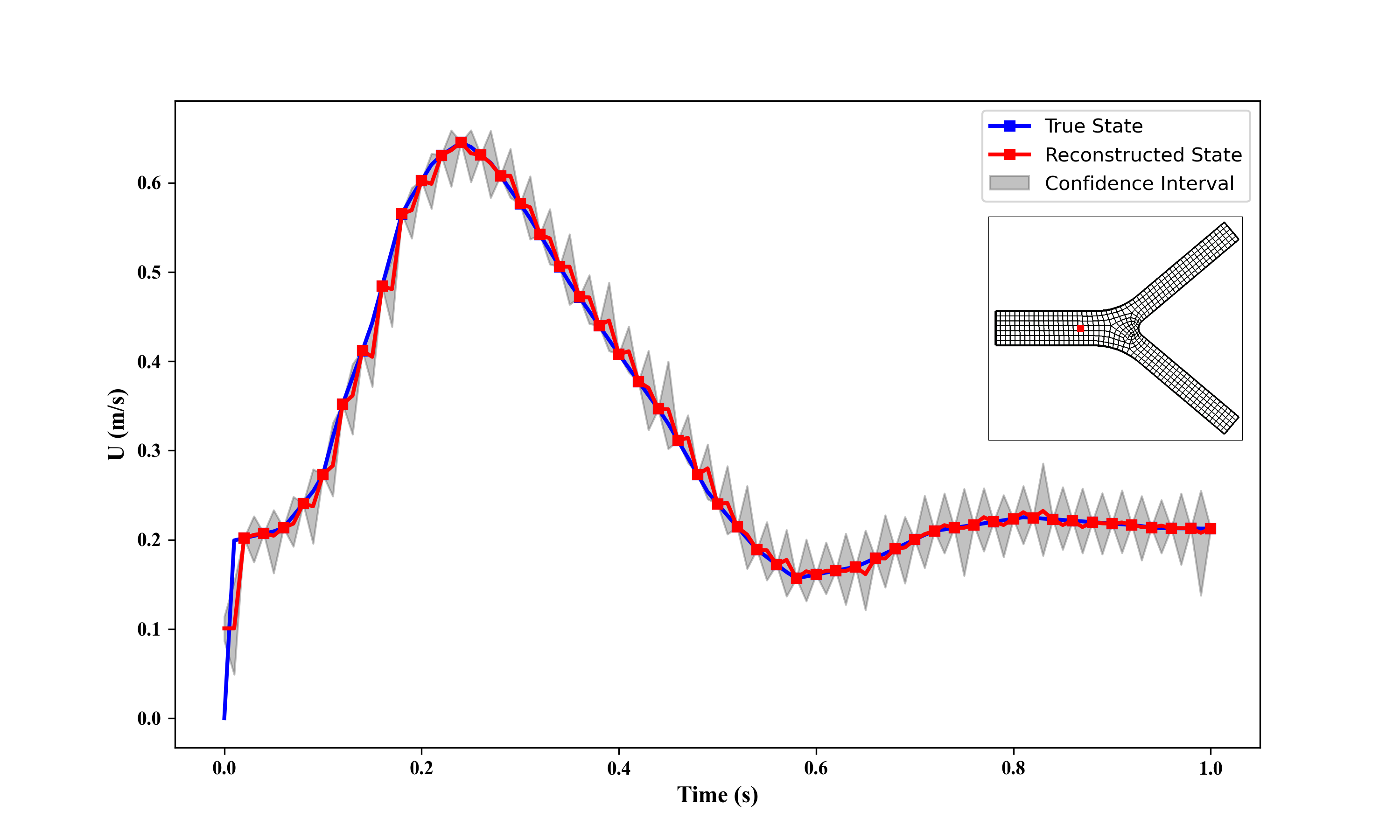}

    \caption{Reconstruction of x-velocity state at the measurement locations ($x = 0.035147$ and $y = 0.007069$)}
    \label{fig:TimeDependentState_timeParameter}
\end{figure}

Fig. \ref{fig:TimeDependentState_timeParameter} compares the reconstructed and true x-velocity state over time at the measurement location (x = 0.035147, y = 0.007069). Similar to the parameter plot, the first two time steps exhibit a high percentage of error due to the absence of measurement data. However, with the incorporation of observations, the reconstructed state rapidly converges to the true state. The confidence interval narrows significantly at observation points, as seen in the constant parameter case. Notably, the states show less discrepancy from the true data at observation intervals, whereas the confidence interval for the parameter does not reach the same level of accuracy, highlighting the added complexity of parameter prediction with EnKF methods.

\begin{figure}[H]
    \centering
    
    \begin{subfigure}{0.7\textwidth}
        \includegraphics[width=\linewidth]{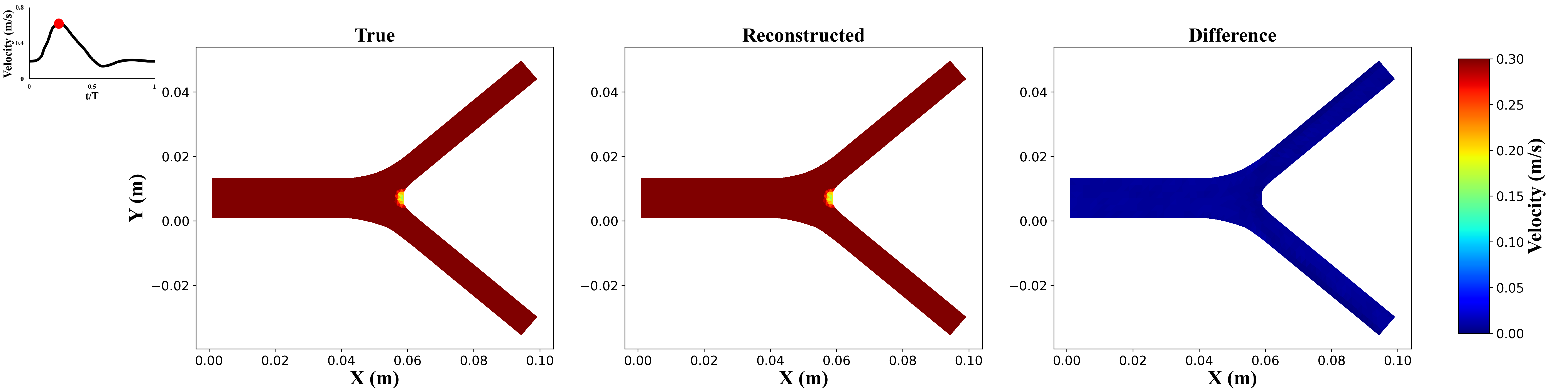}
        \label{fig:Velocity_Contour_diff_timeParameter_at_0.24s}
    \end{subfigure}
    
    \begin{subfigure}{0.7\textwidth}
        \includegraphics[width=\linewidth]{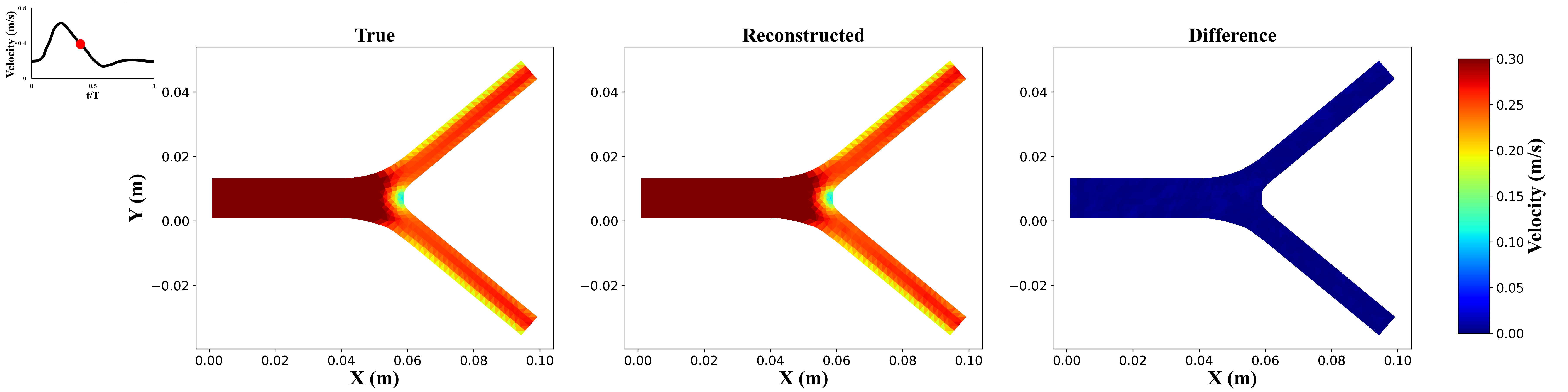} 
        \label{fig:Velocity_Contour_diff_timeParameter_at_0.40s}
    \end{subfigure}
    
    \begin{subfigure}{0.7\textwidth}
        \includegraphics[width=\linewidth]{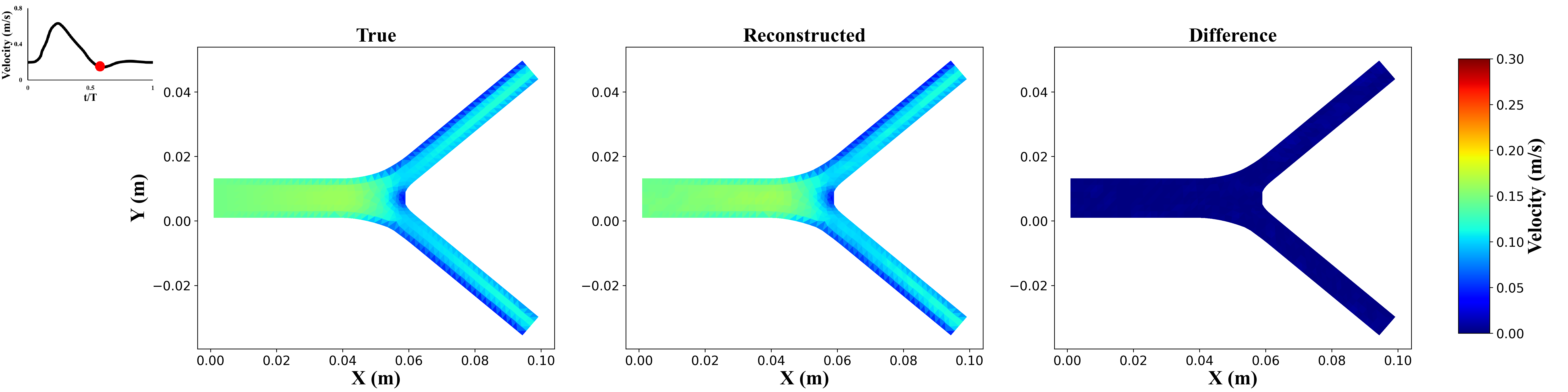} 
        \label{fig:Velocity_Contour_diff_timeParameter_at_0.60s}
    \end{subfigure}
    
    \begin{subfigure}{0.7\textwidth}
        \includegraphics[width=\linewidth]{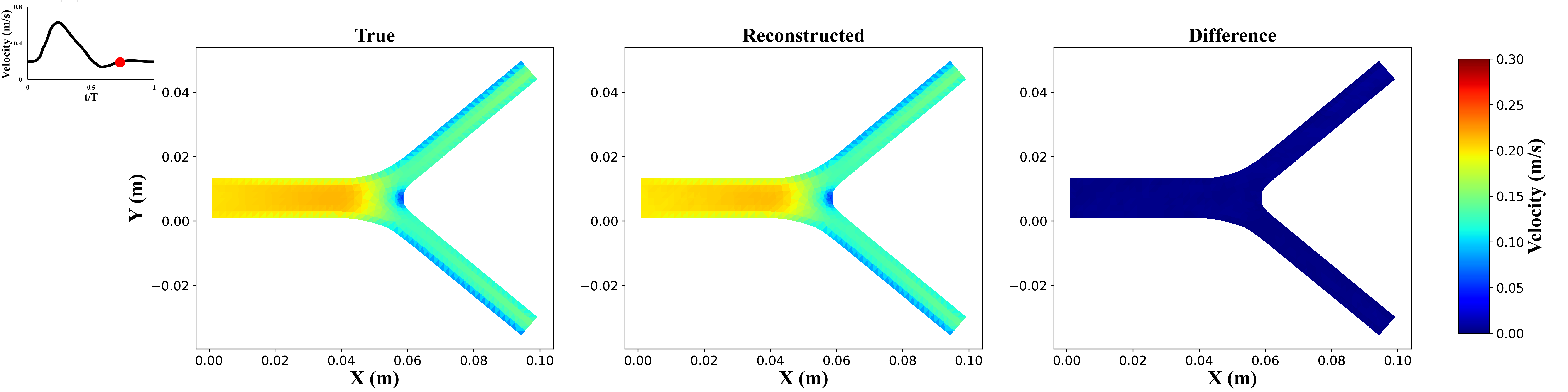} 
        \label{fig:Velocity_Contour_diff_timeParameter_at_0.74s}
    \end{subfigure}    
    
    \begin{subfigure}{0.7\textwidth}
        \includegraphics[width=\linewidth]{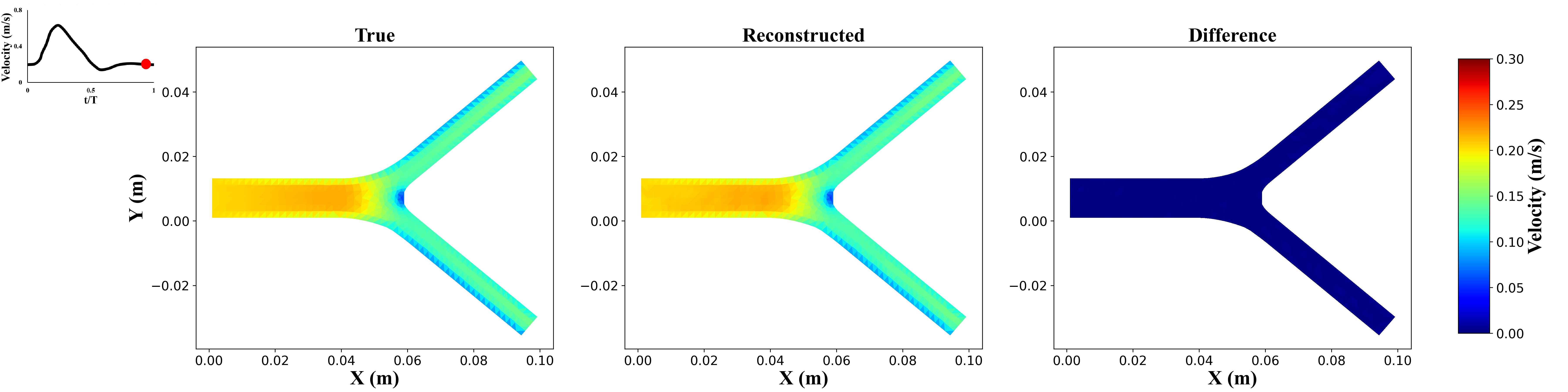} 
        \label{fig:Velocity_Contour_diff_timeParameter_at_0.96s}
    \end{subfigure} 
    
    \caption{Velocity contour comparison between true and reconstructed states and their difference at key cardiac phases for the time-dependent parameter scenario. Each row corresponds to a specific phase, in order: peak systole, maximum deceleration, early diastole, mid diastole, and late diastole.}

    \label{fig:Velocity_Contour_diff_timedependent_comparison}
\end{figure}

Figure \ref{fig:Velocity_Contour_diff_timedependent_comparison} shows velocity magnitude contour plots at five key cardiac phases: peak systole (t/T = 0.24), maximum deceleration (t/T = 0.4), early diastole (t/T = 0.6), mid diastole (t/T = 0.74), and late diastole (t/T = 0.96). Each row displays the contour of the true velocity field (left), the reconstructed field (middle), and their absolute difference (right). Across all time points, the true and reconstructed velocity fields are nearly identical, resulting in minimal differences.

At peak systole, the highest ejection velocity occurs, corresponding to the maximum value of the parameter. Here, high velocity magnitudes are observed throughout the domain, with a noticeable decrease at the bifurcation point. At maximum deceleration, the blood flow sharply decreases following peak systole, marking a shift from high to lower velocities. The velocity in the bifurcations is lower than in the inlet channel, with similar contours for the true and reconstructed fields.

During early diastole, the velocity field significantly decreases, and reverse flow may occur in some regions due to a pressure drop as the heart relaxes. In mid diastole, blood flow stabilizes at a minimal level as the heart chambers refill, with the lowest flow rates typically found in the arteries. Finally, in late diastole, blood flow begins to increase slightly as the heart prepares for the next cycle, with slight acceleration due to atrial contraction. In these last three phases, characterized by low velocities, there are minimal differences between the true and reconstructed fields.

Overall, the EnSISF method demonstrates strong performance in reconstructing the flow field in this time-dependent parameter scenario for the 2D ideal model.

%
%
%
%
\subsubsection{Times-Space-Dependent Parameter}\label{RD:2DTimeSpace}

In the most complex scenario, EnSISF is used to predict a parameter that varies across both time and space. Here, we assume a fully developed boundary profile at the inlet, given by the parabolic formula:
\begin{equation}
\label{eq_ParabVel}
V_{\text{inlet}} = V_{\text{max}} \cdot \left( 1 - \left( \frac{r}{R} \right)^2 \right)
\end{equation}

where \( V_{\text{max}} \) is the maximum velocity in the flow profile, \( r \) is the radial position, and \( R \) is the radius of the inlet boundary. The objective of EnSISF is to predict \( V_{\text{max}} \).

\begin{figure}[H]
    \centering
\includegraphics[width=1\textwidth]{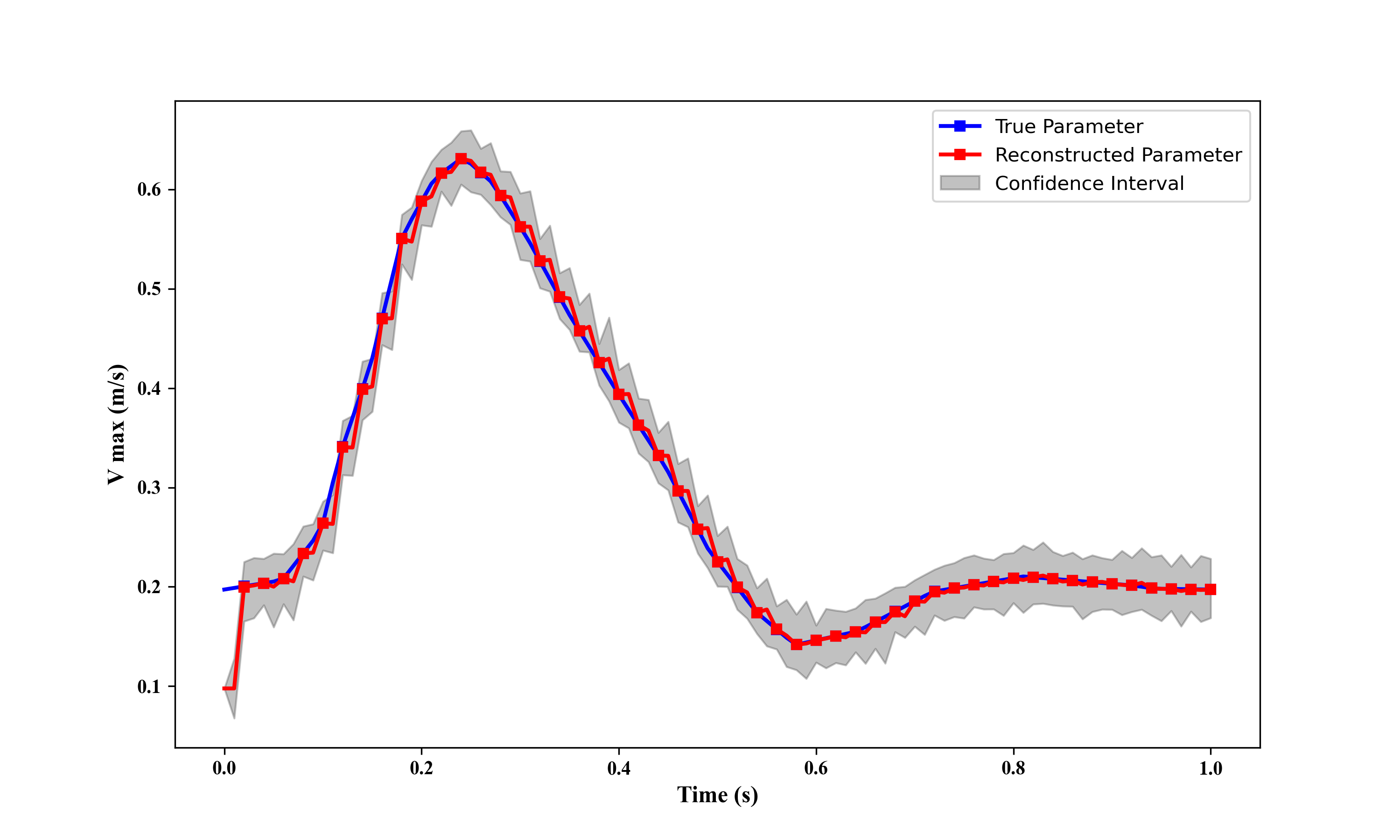}

    \caption{Comparison of true vs. reconstructed time-space-dependent velocity parameter}
    \label{fig:TimeDependentComponentParameter}
\end{figure}

Fig. \ref{fig:TimeDependentComponentParameter} compares the true time-varying maximum velocity with its reconstructed counterpart, along with the confidence interval. Covering a simulation time of 1 second, the plot shows that the initial discrepancy between the true and reconstructed parameters decreases as observations are incorporated from the third time step onward. The confidence interval exhibits similar behavior to the time-dependent case, narrowing consistently at observation points and widening between them. Even in this most complex scenario of our 2D ideal model, the EnSISF achieves precise parameter prediction.

\begin{figure}[H]
    \centering
\includegraphics[width=1\textwidth]{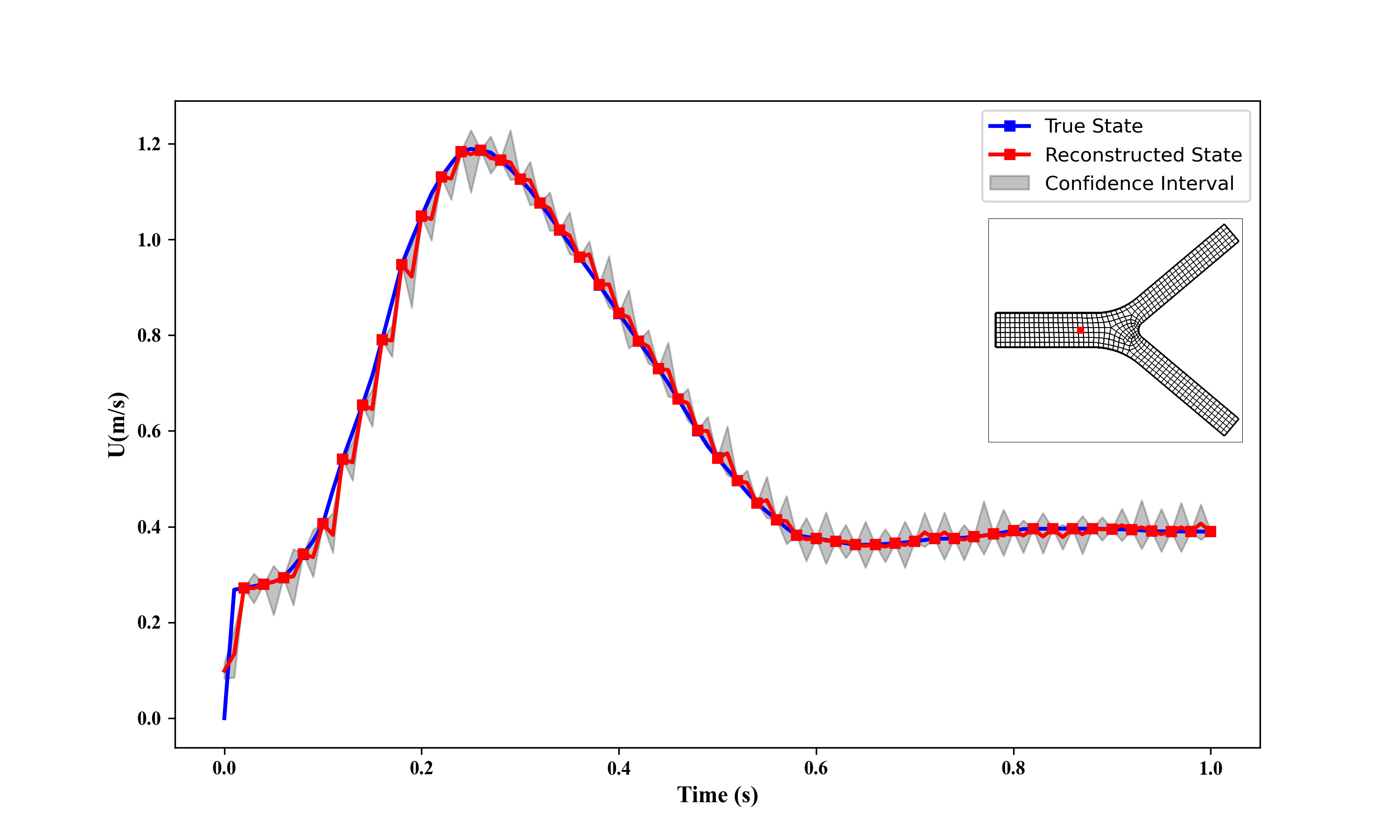}
    \caption{Reconstruction of x-velocity state at the measurement locations ($x = 0.035147$ and $y = 0.007069$)}
    \label{fig:TimeDependentState_timeSpaceParameter}
\end{figure}

Fig. \ref{fig:TimeDependentState_timeSpaceParameter} compares the reconstructed x-velocity with the true x-velocity, as recorded by sensors at the measurement location (x = 0.035147, y = 0.007069) over the simulation time. The plot includes a confidence interval to emphasize the precision and accuracy of the reconstructed variable. The behavior mirrors that of the time-dependent scenario shown in Fig. \ref{fig:TimeDependentState_timeParameter}, and a similar explanation applies to this figure.

\begin{figure}[H]
    \centering
    
    \begin{subfigure}{0.7\textwidth}
        \includegraphics[width=\linewidth, height=0.25\textheight, keepaspectratio]{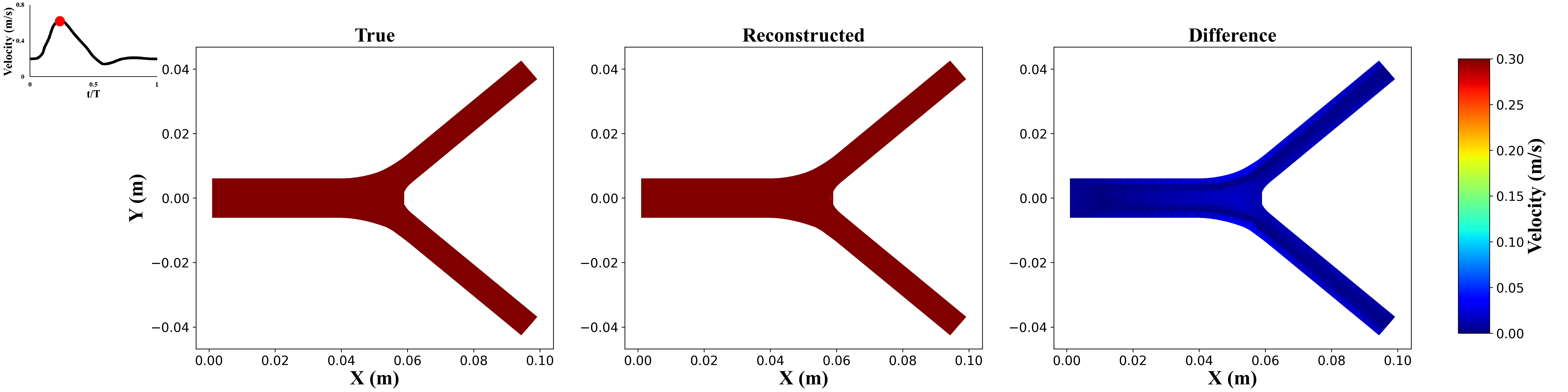}
        \label{fig:Velocity_Contour_diff_at_0.24s}
    \end{subfigure}
    
    \begin{subfigure}{0.7\textwidth}
        \includegraphics[width=\linewidth, height=0.25\textheight, keepaspectratio]{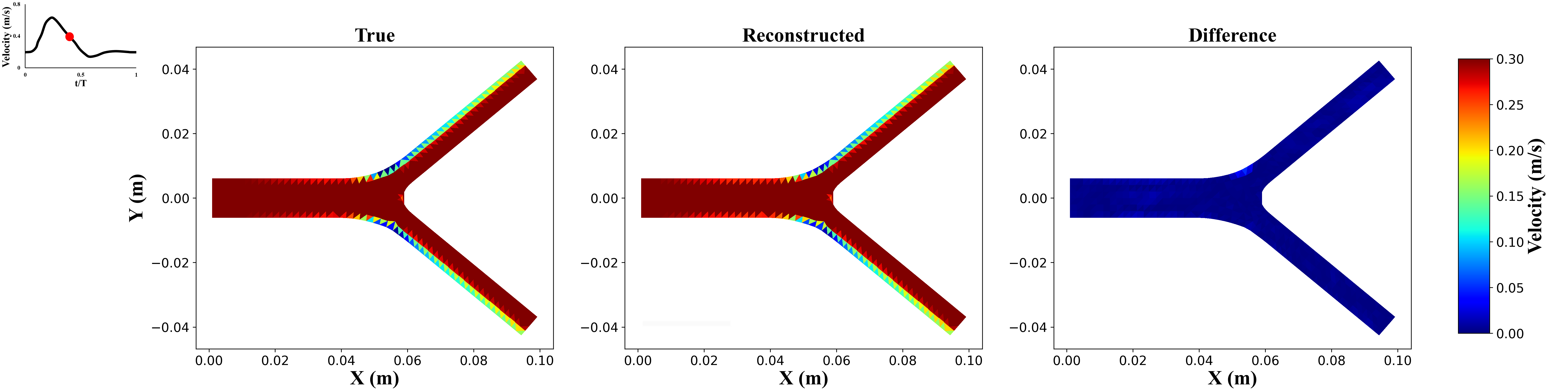} 
        \label{fig:Velocity_Contour_diff_at_0.40s}
    \end{subfigure}
    
    \begin{subfigure}{0.7\textwidth}
        \includegraphics[width=\linewidth, height=0.25\textheight, keepaspectratio]{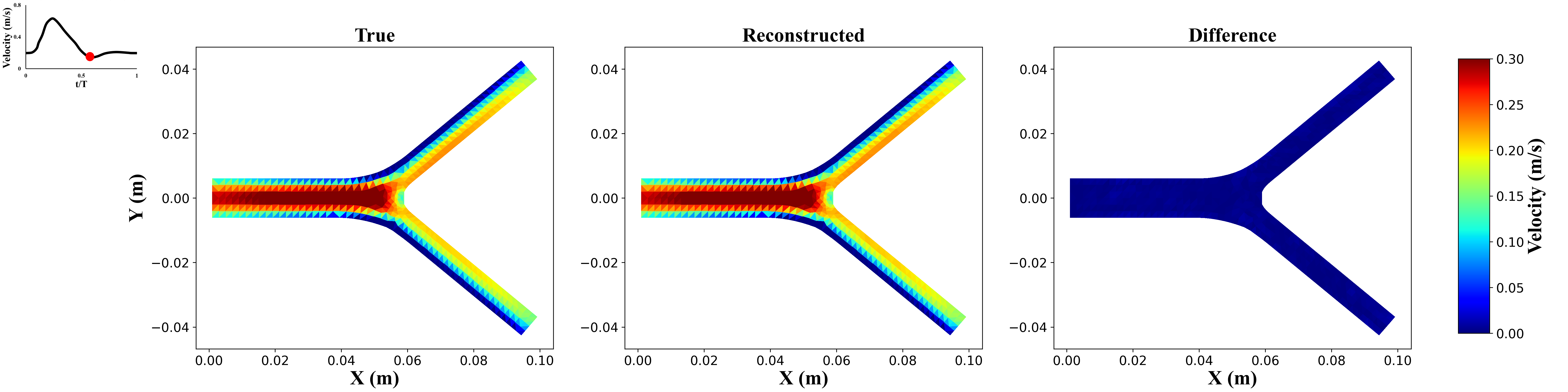} 
        \label{fig:Velocity_Contour_diff_at_0.60s}
    \end{subfigure}
    
    \begin{subfigure}{0.7\textwidth}
        \includegraphics[width=\linewidth, height=0.25\textheight, keepaspectratio]{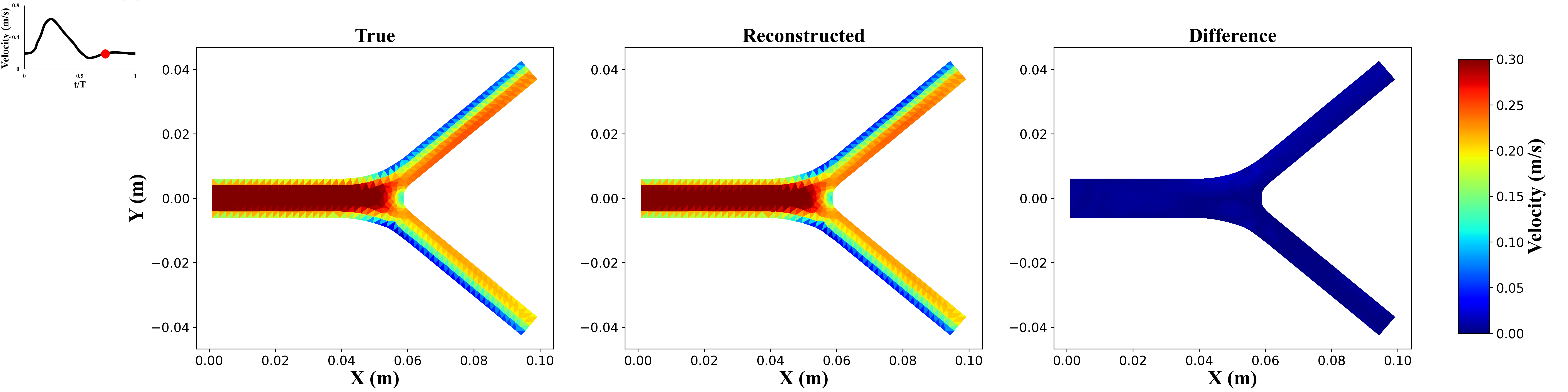} 
        \label{fig:Velocity_Contour_diff_at_0.74s}
    \end{subfigure}    

    \begin{subfigure}{0.7\textwidth}
        \includegraphics[width=\linewidth, height=0.25\textheight, keepaspectratio]{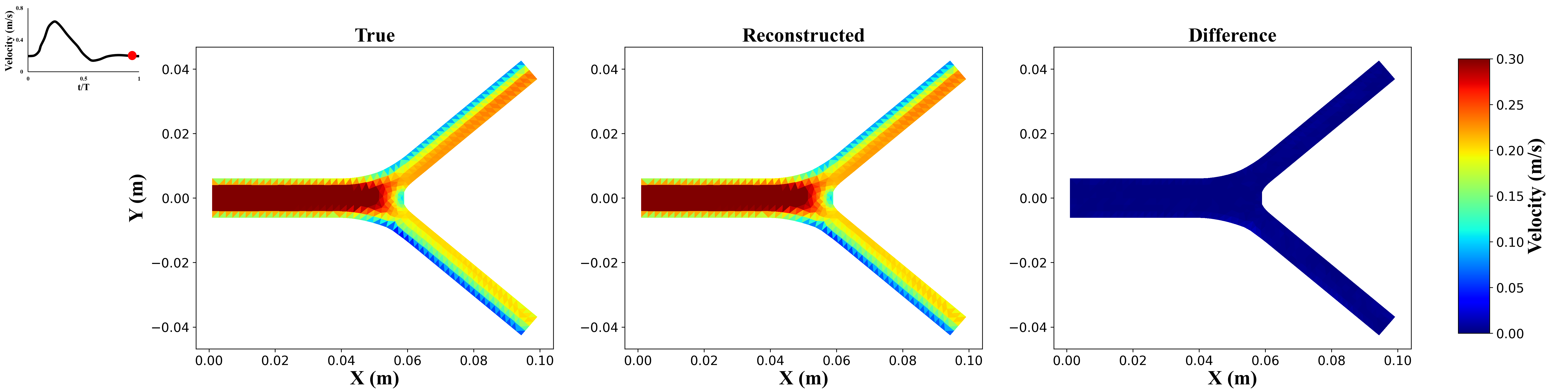} 
        \label{fig:Velocity_Contour_diff_at_0.96s}
    \end{subfigure}  
    
    \caption{Velocity contour comparison between true and reconstructed states and their difference at key cardiac phases for the time-space-dependent parameter scenario. Each row corresponds to a specific phase, in order: peak systole, maximum deceleration, early diastole, mid diastole, and late diastole.}
    \label{fig:Velocity_Contour_timeSpaceDependent_diff_comparison}
\end{figure}

Similar to the previous scenario, Fig. \ref{fig:Velocity_Contour_timeSpaceDependent_diff_comparison} presents the true and reconstructed velocity fields for each cardiac phase in separate subplots. Here, a parabolic inlet velocity profile results in a non-uniform velocity distribution in the inlet channel, unlike the time-dependent scenario. Despite this, the velocity difference between the fields remains minimal, generally less than 0.05 m/s at those time points. However, in contrast to the previous scenario, peak systole shows a higher error, particularly near the wall boundaries and bifurcation point. This increased error arises from the complexity of this time-space-dependent scenario, the presence of high non-linearity near the wall regions, and the fact that peak systole corresponds to a fully turbulent field, whereas the EnSISF model uses a laminar solver to manage computational demands.

The error level decreases at maximum deceleration but is still somewhat higher than in the time-dependent parameter scenario. These errors tend to reduce further in cardiac phases associated with lower velocity magnitudes. Despite these complexities, the model effectively reconstructs the flow field with acceptable accuracy.

\begin{figure}[H]
     \centering

     \begin{subfigure}{0.45\textwidth}
         \includegraphics[width=\linewidth]{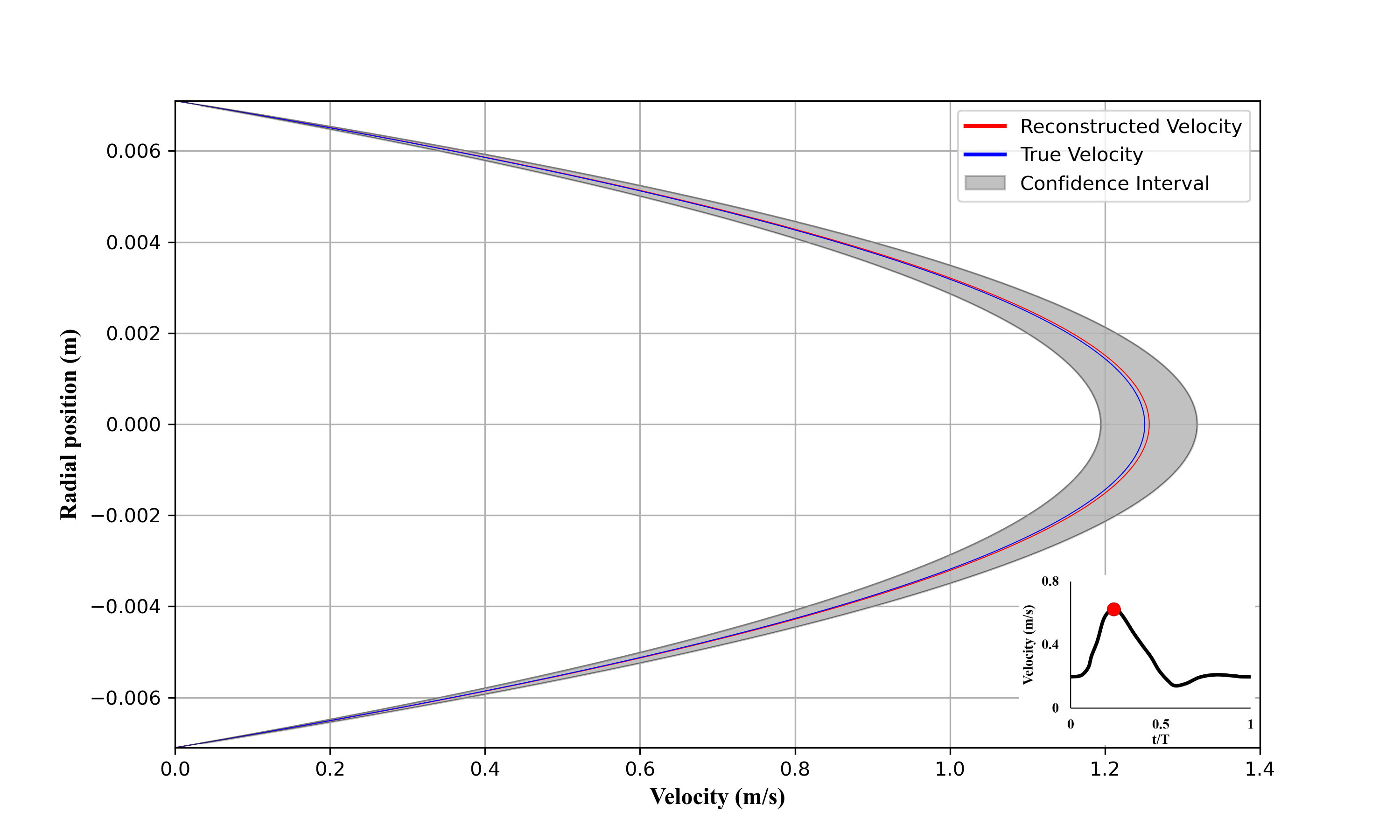}
         \label{fig:VelocityProfile_t_0_25s}
     \end{subfigure}
     \hfill
     \begin{subfigure}{0.45\textwidth}
         \includegraphics[width=\linewidth]{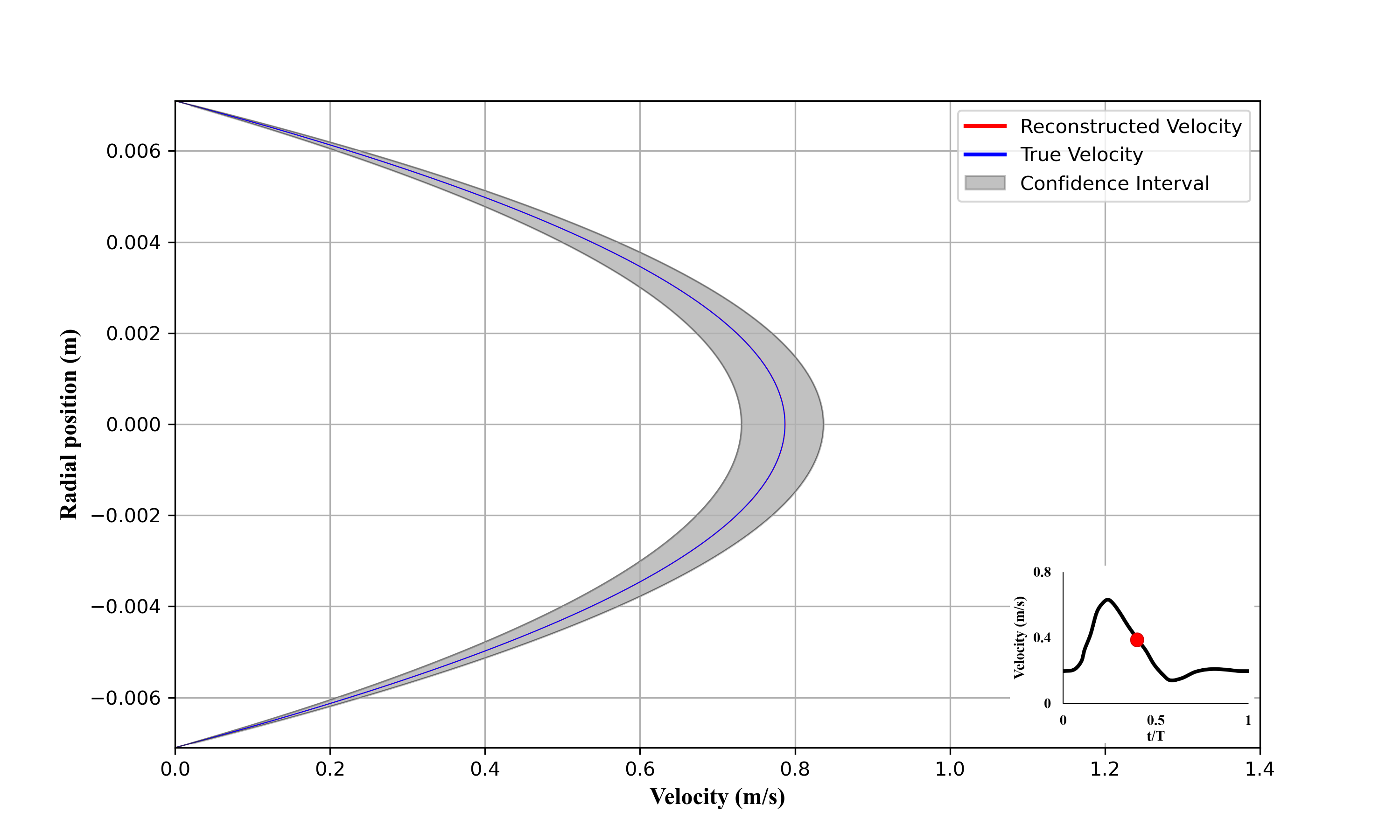} 
         \label{fig:VelocityProfile_t_0_40s}
     \end{subfigure}

     \vspace{0.3cm}

     \begin{subfigure}{0.45\textwidth}
         \includegraphics[width=\linewidth]{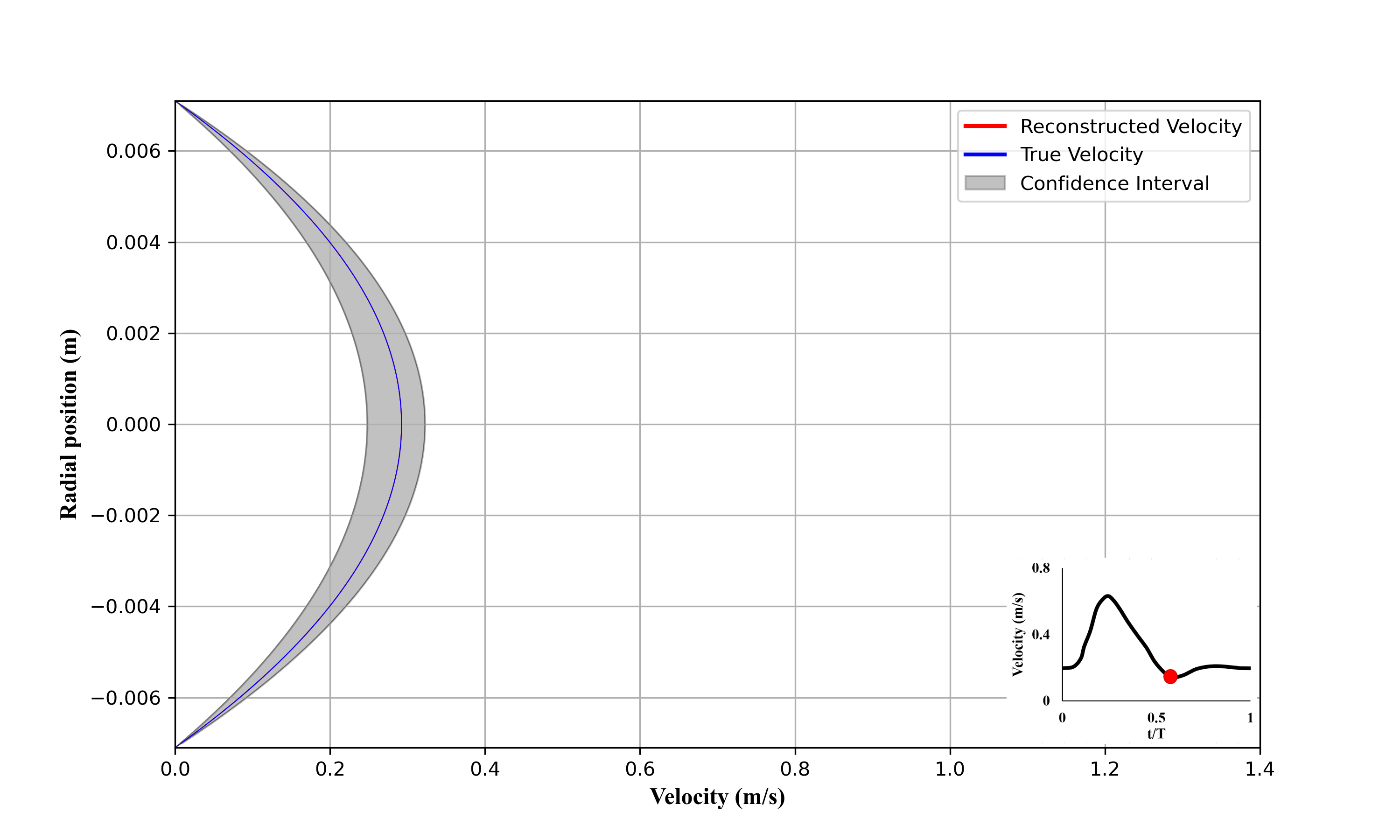} 
         \label{fig:VelocityProfile_t_0_60s}
     \end{subfigure}
     \hfill
     \begin{subfigure}{0.45\textwidth}
         \includegraphics[width=\linewidth]{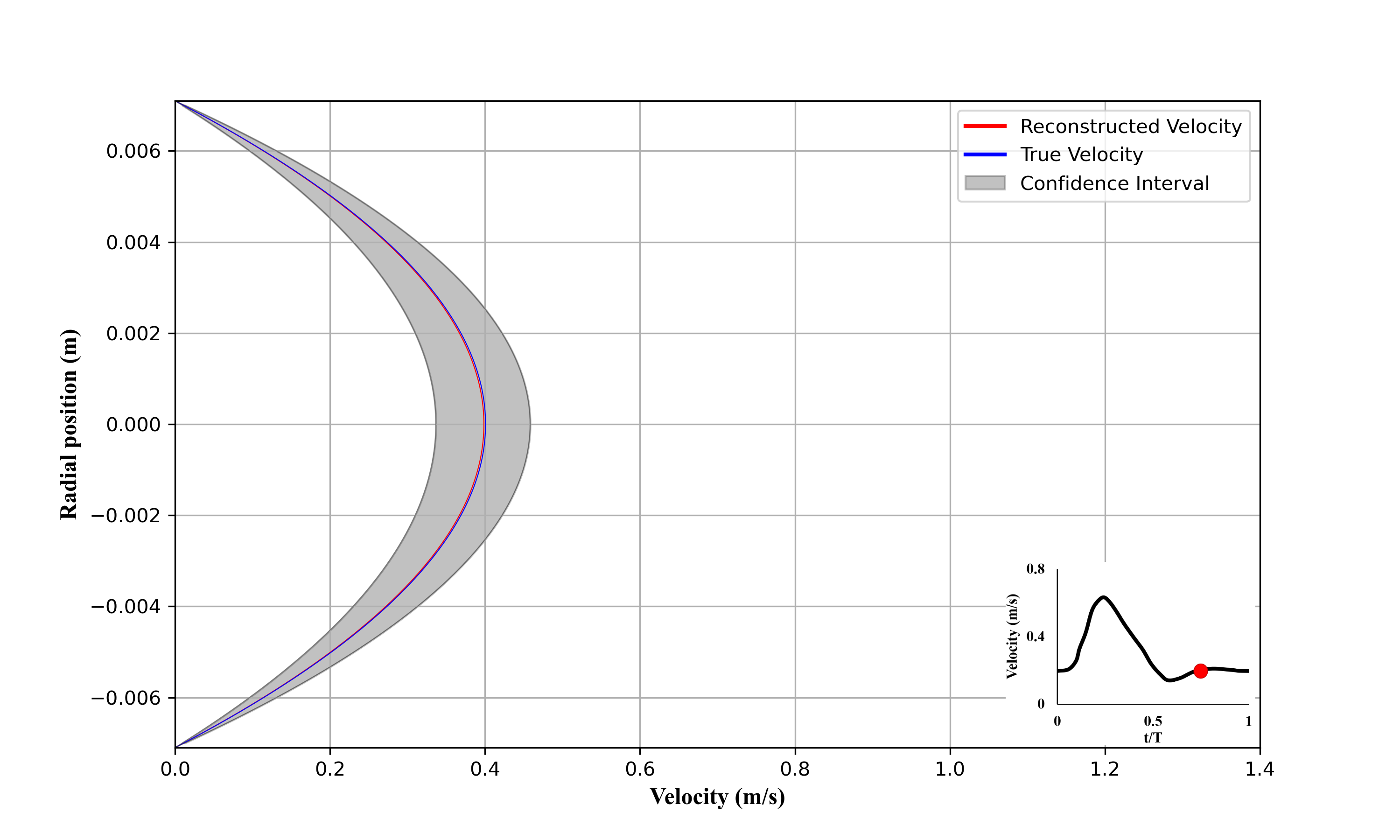} 
         \label{fig:VelocityProfile_t_0_75s}
     \end{subfigure}

     \vspace{0.3cm}

     \begin{subfigure}{0.45\textwidth}
         \includegraphics[width=\linewidth]{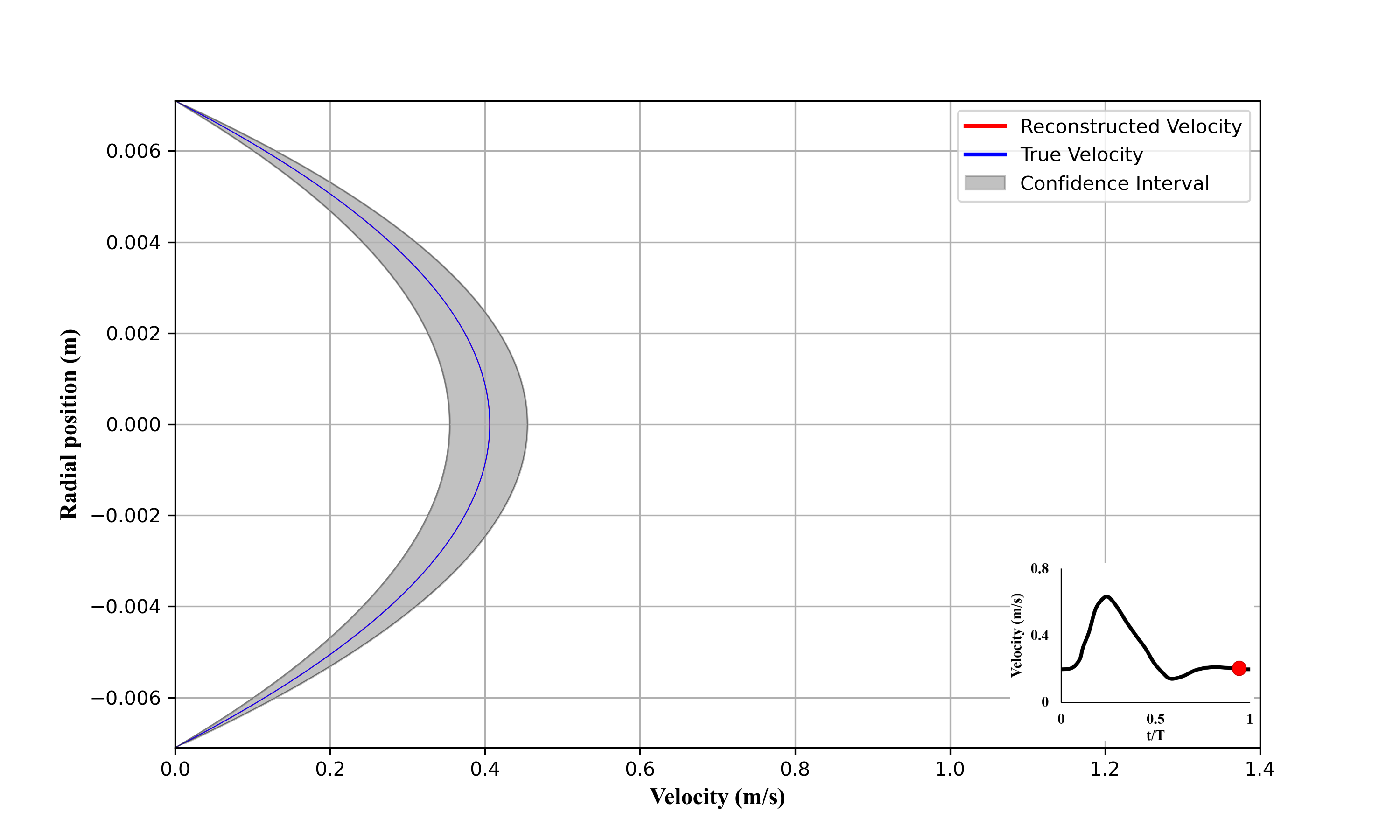} 
         \label{fig:VelocityProfile_t_0_90s}
     \end{subfigure}

     \caption{Inlet velocity profile comparison between true and reconstructed states and their difference at key cardiac phases for the time-space-dependent parameter scenario. Each image corresponds to a specific phase, in order: peak systole, maximum deceleration, early diastole, mid diastole, and late diastole.}
     \label{fig:velocityProfile_comparison}
\end{figure}

Fig. \ref{fig:velocityProfile_comparison} shows the inlet velocity profiles at specified time points for both true and reconstructed velocities, along with the confidence interval. The EnSISF accurately predicts the parameter with an acceptable error margin across these intervals. However, at peak systole, the time of highest velocity, there is a slight discrepancy between the true and reconstructed parameter. Additionally, the confidence interval at this point covers a wider range, whereas for the other time points, it is slightly narrower. Later time instances show improved parameter estimation, aligning well with the field data presented in Figure \ref{fig:Velocity_Contour_timeSpaceDependent_diff_comparison}.

    
    
    
    

    

\subsection{3D Patient-Specific Model}\label{3DTimeSpace}

In our final assessment, we applied EnSISF to a 3D patient-specific model, investigating the most complex case of a time-space-dependent parameter. A relative error of 7.37\% was observed with an observation interval of 0.02 seconds.

\begin{figure}[H]
    \centering
    \includegraphics[width=\linewidth, height=0.8\textheight, keepaspectratio]{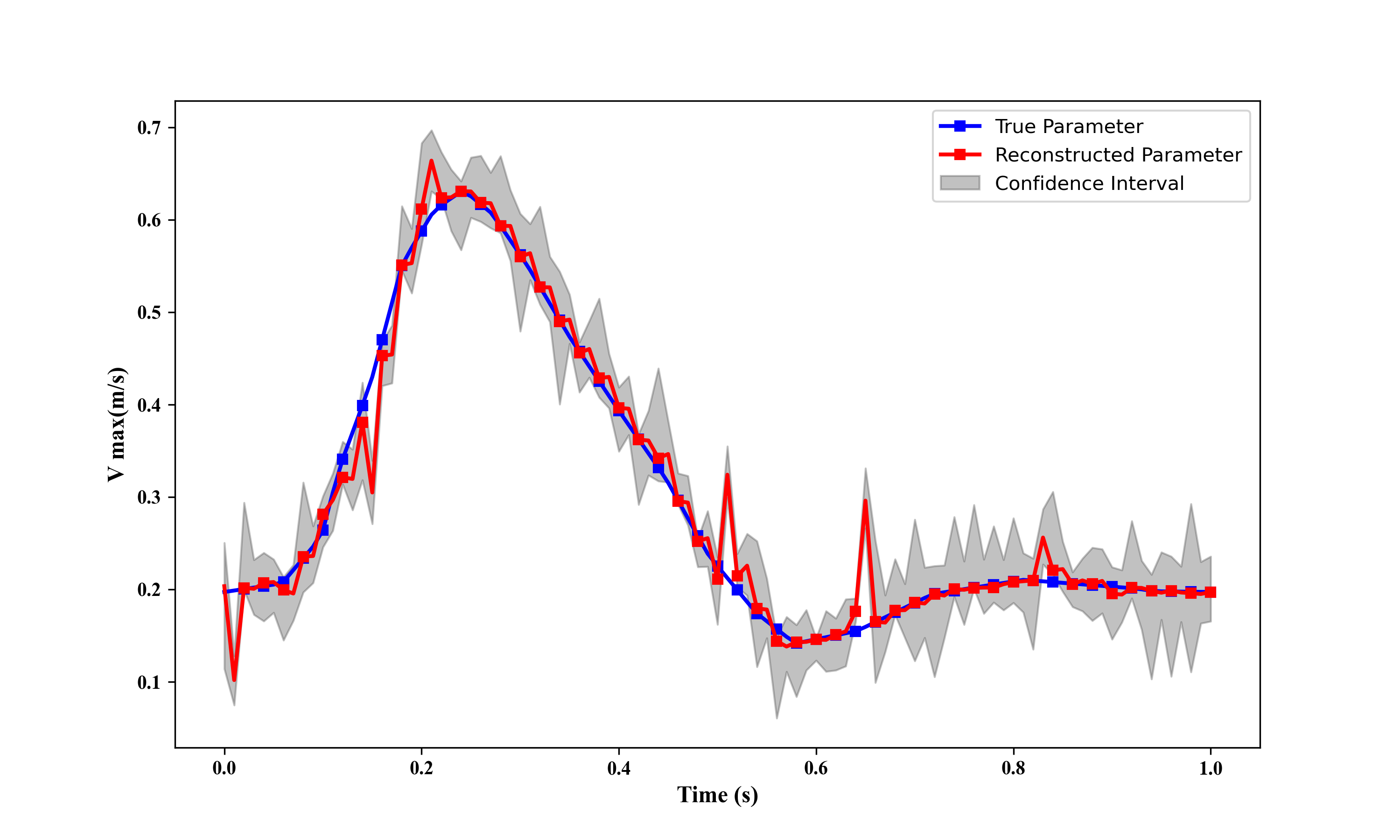}
    \caption{Comparison of true vs. reconstructed time-space-dependent parameter 3D patient-specific abdominal aorta}
    \label{TrueReconstructedParameter3D}
\end{figure}

Fig. \ref{TrueReconstructedParameter3D} shows a comparison between the true and reconstructed values of a time-varying boundary condition parameter over the cardiac cycle. In this complex case, a greater degree of deviation was noted. At two specific time points, t = 0.65 sec and 0.83 sec (intervals without observations), significant divergence in the predicted parameter occurred. However, at later time steps, with measurement assimilation, the model converged to the actual parameter values. This case exhibited higher oscillations than the 2D ideal model, with larger confidence intervals that did not always coincide with the true values at certain time points, yet still demonstrated relatively good accuracy in parameter prediction.

\begin{figure}[H]
    \centering
    \includegraphics[width=\linewidth, height=0.8\textheight, keepaspectratio]{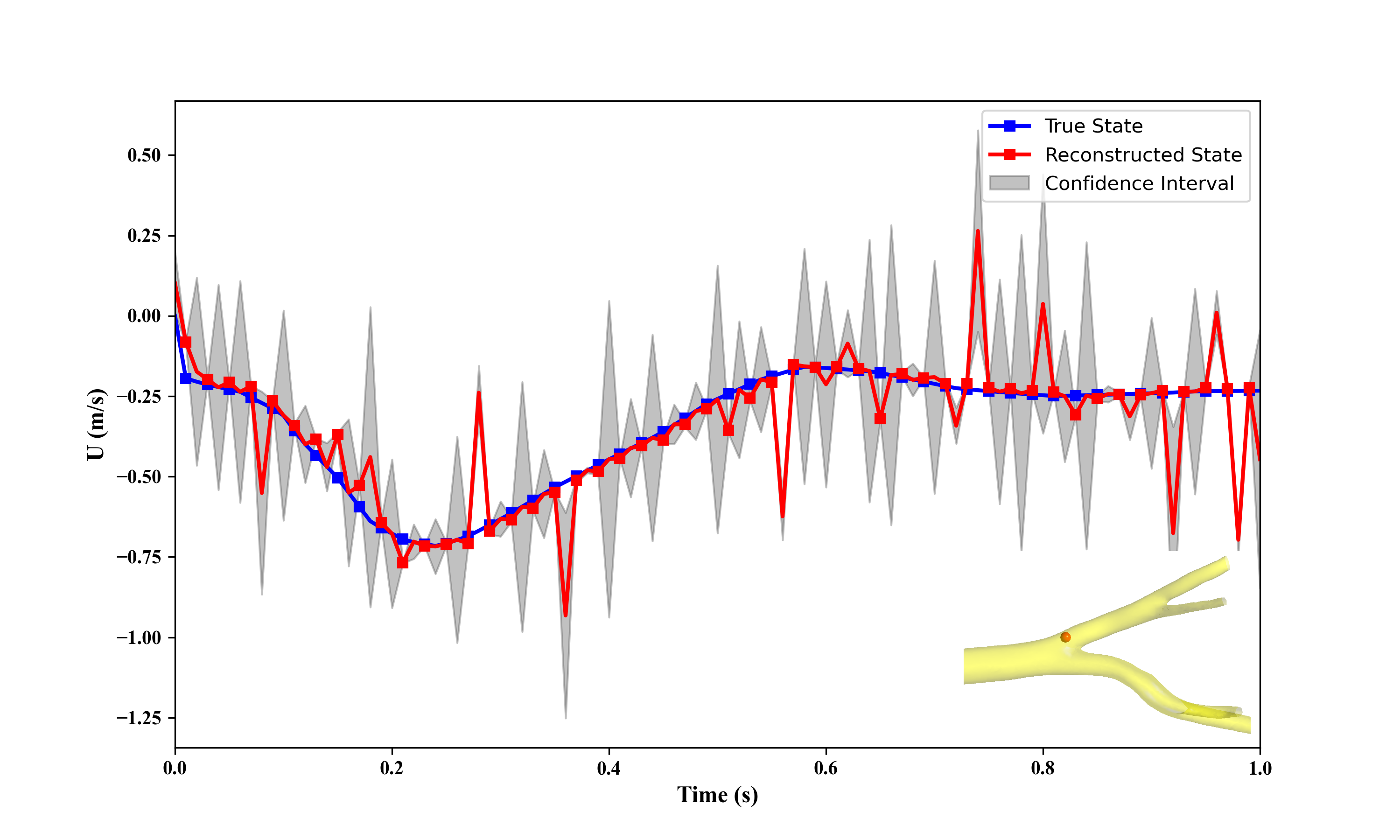}
    \caption{Reconstruction of z-velocity state at the measurement locations (x = -0.0117, y = 0.0031 and z = -0.0397) for the 3D patient-specific case}
    \label{TrueReconstructedState3D}
\end{figure}

In Fig. \ref{TrueReconstructedState3D}, a comparison of the true and reconstructed z-velocity over time is presented based on sensor data from the location (x = -0.0117, y = 0.0031, z = -0.0397). The deviation between true and reconstructed values is relatively high in intervals without observations; however, with data assimilation, the model aligns closely with the true values during observation intervals. This behavior reflects the complexity of real-world models, which require a greater number of spatial measurement points."



\begin{figure}[H]
    \centering
    \includegraphics[width=\linewidth, height=0.8\textheight, keepaspectratio]{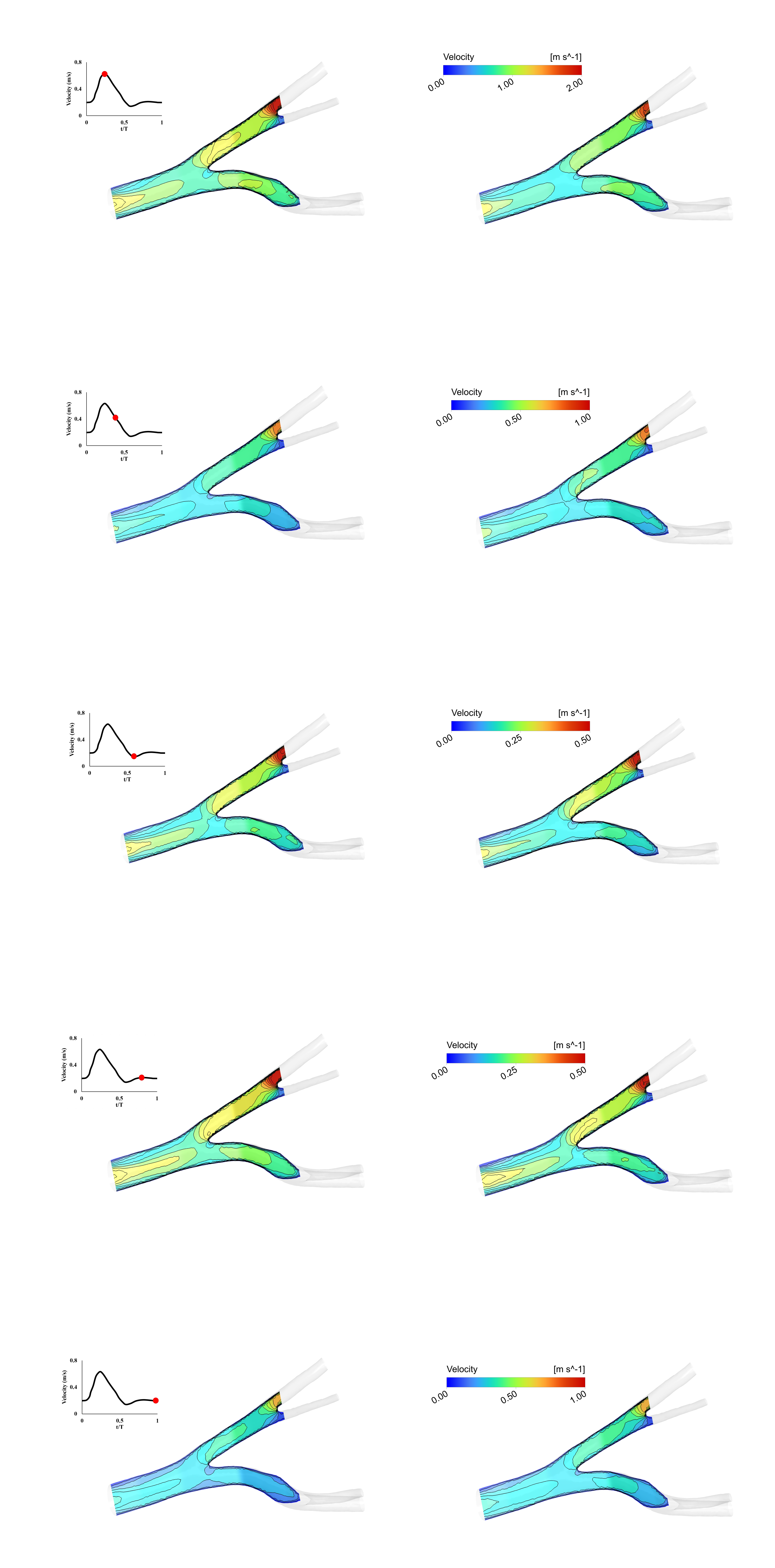}
    \caption{Comparison of velocity contours between pseudo-experimental data (left) and the predicted mean distribution (right) at key cardiac phases: peak systole, maximum deceleration, early diastole, mid diastole, and late diastole}
    \label{comparison_velocity}
\end{figure}

\begin{figure}[H]
    \centering
    \includegraphics[width=\linewidth, height=0.8\textheight, keepaspectratio]{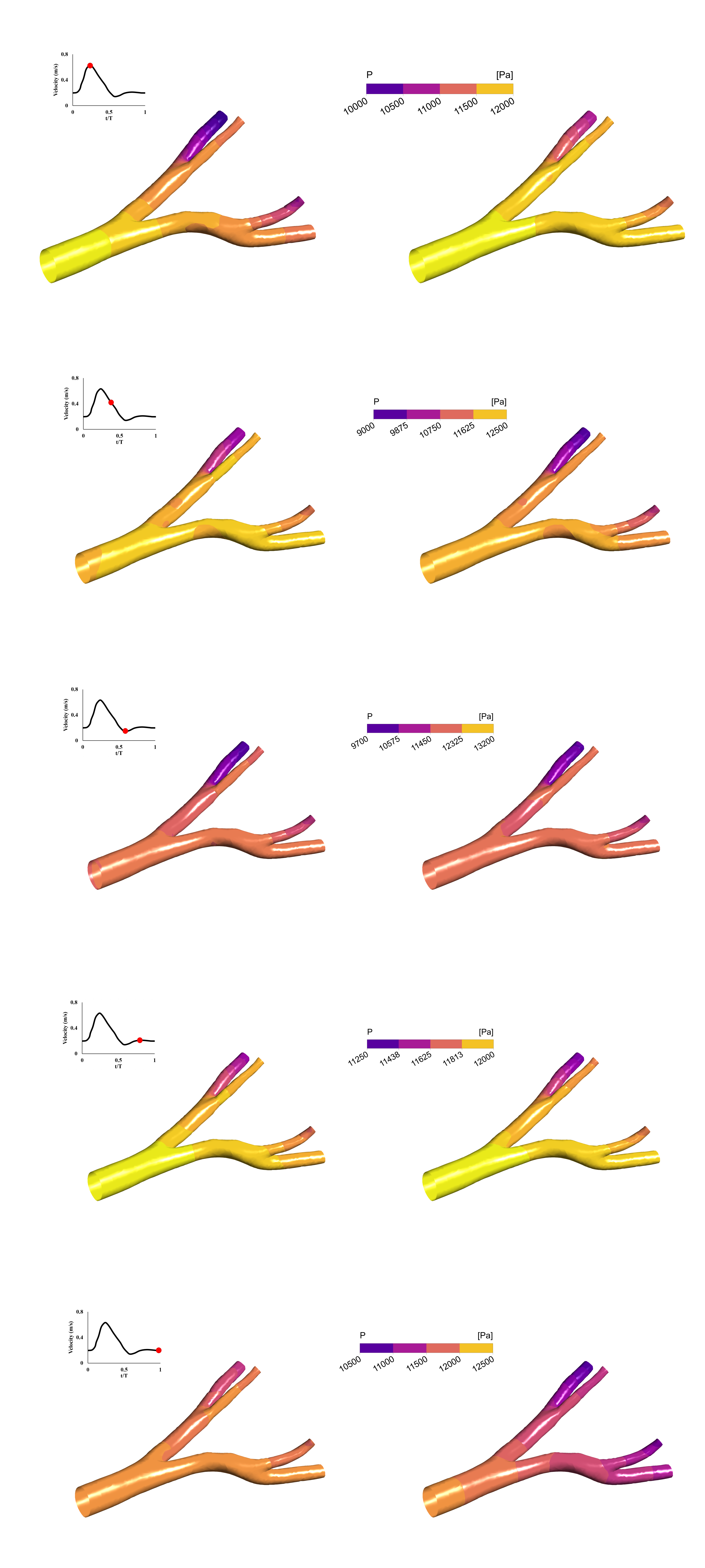}
    \caption{Comparison of pressure contours between pseudo-experimental data (left) and the predicted mean distribution (right) at key cardiac phases: peak systole, maximum deceleration, early diastole, mid diastole, and late diastole}
    \label{comparison_pressure}
\end{figure}

Figs. \ref{comparison_velocity} and \ref{comparison_pressure} compare the velocity and pressure contours between pseudo-experimental (left column) and reconstructed fields (right column) at key cardiac intervals. At peak systole, corresponding to the highest velocity interval, there is a noticeable discrepancy between the fields, particularly in the inlet channel and at the bifurcations. This is also reflected in the pressure field at the same time point. Similar to the 2D results in Fig. \ref{fig:Velocity_Contour_timeSpaceDependent_diff_comparison}, this highlights the sensitivity of the EnSISF method to velocity magnitude, with greater deviation observed at higher velocities. Additionally, it’s important to note that turbulent flow occurs at this time interval, whereas the forward solver uses a laminar model to manage computational demands.

Despite this, the overall comparability of results and the EnSISF's ability to capture both trends and values for the velocity and pressure fields indicate that a high-fidelity solver may not be necessary for the forward solver in the EnKF method. Discrepancies decrease at other time points. During maximum deceleration, slight differences remain, especially in the bifurcation regions where the pressure field shows minor variation. Throughout diastole, the velocity contours align closely, as does the pressure field, demonstrating strong agreement overall.

\section{Conclusion}\label{sec:conc}

This study aimed to utilize an advanced ensemble Kalman filter variant, EnSISF, to estimate unknown boundary velocity profiles in cardiovascular models. Recognizing the limitations of in-vivo measurement data, such as 4D flow MRI, which often suffer from low resolution and noise, our approach not only addresses these uncertainties but also accurately predicts unknown profiles for clinical applications. Accurately determining the velocity profile enables further calculation of wall shear indices, which are critical for diagnosing and treating cardiovascular diseases like aneurysms and atherosclerosis.

Simulations were conducted for constant, time-dependent, and time-space-dependent parameters using a 2D ideal model of the abdominal aorta, eventually extending to predict a time-space-dependent parameter in a 3D patient-specific model. In our 2D simulations, we observed relative errors as low as 0.996\% for constant boundary conditions over a 0.02-second observation span, which increased to 2.63\% for time-dependent boundaries and 2.61\% for time-space-dependent boundaries, indicating excellent precision in boundary estimation for 2D cases. In the 3D patient-specific model, where both time and space dependencies were considered, the relative error reached 7.37\% over a similar observation span, highlighting the challenges in 3D cardiovascular simulations.

The level of discrepancy between the true and reconstructed fields varied across specific time points, with peak systole showing the highest discrepancy due to the highest inlet velocity—reflecting the use of a simplified laminar solver for the forward update of EnSISF. Other time points, associated with lower velocities, exhibited minimal differences in both the field and predicted inlet velocity. These findings validate the accuracy of EnSISF in reconstructing the velocity field over time, underscoring its potential for patient-specific hemodynamic applications.
\section*{Acknowledgements}

Professor Giovanni Stabile acknowledges the financial support under the National Recovery and Resilience Plan (NRRP), Mission 4, Component 2, Investment 1.1, Call for tender No. 1409 published on 14.9.2022 by the Italian Ministry of University and Research (MUR), funded by the European Union – NextGenerationEU– Project Title ROMEU – CUP P2022FEZS3 - Grant Assignment Decree No. 1379 adopted on 01/09/2023 by the Italian Ministry of Ministry of University and Research (MUR). We acknowledge the support provided by PRIN “FaReX - Full and Reduced order modelling of coupled systems: focus on non-matching methods and automatic learning” project, PNRR NGE iNEST “Interconnected Nord-Est Innovation Ecosystem” project, INdAM-GNCS 2019–2020 projects and PON “Research and Innovation on Green related issues” FSE REACT-EU 2021 project.
\section*{Data Availability}
All data generated or analyzed during this study are included in this published article.
\bibliographystyle{unsrt}
\bibliography{biblio.bib}
\end{document}